\documentclass[11pt]{amsart}

\usepackage{amsmath,amsthm,amsfonts,amssymb,mathtools}
\usepackage{dsfont}
\usepackage{enumerate}
\usepackage{xcolor}
\usepackage{tikz}
\usepackage[colorlinks=true, linkcolor=blue, citecolor=magenta, menucolor=black]{hyperref}

\makeatletter
\@ifundefined{subjclassname@2020}{%
    \@namedef{subjclassname@2020}{\textup{2020} Mathematics Subject Classification}
}{}
\makeatother

\numberwithin{equation}{section}
\makeatletter
\def\subsection{\@startsection{subsection}{2}%
  \z@{.7\linespacing\@plus.7\linespacing}{.3\linespacing}%
  {\normalfont\bfseries}}
\makeatother

\renewcommand{\phi}{\varphi}

\newcommand{\zz}{\mathbb Z}
\newcommand{\nn}{\mathbb N}
\newcommand{\rr}{\mathbb R}
\newcommand{\qq}{\mathbb Q}
\newcommand{\cc}{\mathbb C}

\newcommand{\un}{\mathds 1}

\newcommand{\A}{\mathcal A}
\renewcommand{\H}{\mathcal H}
\newcommand{\U}{\mathcal U}
\newcommand{\Z}{\mathcal Z}

\newcommand{\pt}{\otimes}
\newcommand{\ovt}{\mathbin{\overline{\otimes}}}

\DeclareMathOperator{\LL}{L}
\DeclareMathOperator{\I}{I}
\DeclareMathOperator{\II}{II}
\DeclareMathOperator{\Aut}{Aut}
\DeclareMathOperator{\supp}{supp}
\DeclareMathOperator{\sgn}{sgn}
\DeclareMathOperator{\Tr}{Tr}
\DeclareMathOperator{\rk}{rk}
\DeclareMathOperator{\rank}{rank}
\DeclareMathOperator{\lcm}{lcm}

\newcommand{\abs}[1]{\left\lvert #1\right\rvert}
\newcommand{\norm}[1]{\left\Vert #1\right\Vert}
\newcommand{\ps}[2]{\left\langle #1,#2\right\rangle}
\newcommand{\ovl}[1]{\overline{#1}}

\newcommand{\Wstar}{\operatorname{W}^*}
\newcommand{\Cstar}{\operatorname{C}^*}
\newcommand{\Sym}{\mathfrak S}
\newcommand{\actson}{\curvearrowright}

\newcommand{\rB}{\operatorname{B}}
\newcommand{\rC}{\operatorname{C}}

\newcommand{\rL}{\operatorname{L}}
\newcommand{\rE}{\operatorname{E}}
\newcommand{\rM}{\operatorname{M}}
\newcommand{\Hk}[2]{\rC_c(#2\backslash #1/#2)}

\makeatletter
\newcommand{\@optsub}[1]{\if\relax\detokenize{#1}\relax\else_{#1}\fi}
\newcommand{\csupp}[2][]{c\@optsub{#1}(#2)}
\newcommand{\zat}[1][]{z^{\mathrm{a}}\@optsub{#1}}
\newcommand{\zdiff}[1][]{z^{\mathrm{d}}\@optsub{#1}}
\makeatother

\newtheorem{The}{Theorem}[section]
\newtheorem{Prop}[The]{Proposition}
\newtheorem{Lem}[The]{Lemma}
\newtheorem{fact}[The]{Fact}
\newtheorem{claims}[The]{Claims}
\newtheorem{Cor}[The]{Corollary}
\newtheorem{mainthm}{Theorem}

\theoremstyle{definition}
\newtheorem{Def}[The]{Definition}

\newtheorem{Ex}[The]{Example}

\makeatletter
\def\@listI{\leftmargin\leftmargini \parsep\z@skip
  \topsep\listisep \itemsep 2\p@\relax}
\let\@listi\@listI
\makeatother

\begin{document}
\title[The type \textup{I} conjecture for groups acting on trees]{The type \textup{I} conjecture for locally compact groups acting on trees}
\author{Cyril Houdayer}
\author{Basile Morando}

\address{\'Ecole normale sup\'erieure \\ D\'epartement de math\'ematiques et applications \\ Universit\'e Paris-Saclay \\ 45 rue d'Ulm \\ 75230 Paris Cedex 05 \\ France}
\email{cyril.houdayer@ens.psl.eu}

\address{\'Ecole normale sup\'erieure \\ D\'epartement de math\'ematiques et applications \\ 45 rue d'Ulm \\ 75230 Paris Cedex 05 \\ France}
\email{basile.morando@ens.psl.eu}

\thanks{Research supported by ERC Advanced Grant NET 101141693.}


\begin{abstract}
    We prove that a closed group of automorphisms acting minimally on a locally finite tree and two-transitively on its boundary is uniformly admissible, hence of type~$\I$. This solves Nebbia's conjecture in a strong form. Together with the work of Caprace, Kalantar and Monod, this shows that for every closed nonamenable subgroup $G$ of $\Aut(T)$ acting minimally on a locally finite tree $T$, $G$ is of type~$\I$ if and only if the action on $\partial T$ is two-transitive. The key ingredient is the following operator algebraic result: if the projections of two atomic von Neumann algebras $A_1,A_2$ have integer traces, then every amenable amalgamated free product $A_1\ast_BA_2$ with $B$ atomic is of type~$\I$.
\end{abstract}

\maketitle
\vspace{2\baselineskip}
\raggedbottom

\section{Introduction and statement of the main results}
\label{sec:introduction}

The type~$\I$ property goes back to Murray and von Neumann, whose theory of rings of operators (today's von Neumann algebras) grew in part out of their ``\textit{attempts to generalise the theory of unitary group-representations essentially beyond their classical frame}'' \cite{MvN36}. This classical frame, namely Fourier analysis on abelian groups and the Peter--Weyl theory of compact groups, is one in which every representation is uniquely built from irreducible ones. Type~$\I$ groups are precisely those for which it survives (see Theorem~\ref{thm:glimm} for equivalent definitions). Determining whether a given group is of type~$\I$ is therefore, in the words of Auslander and Moore, ``\textit{one of the fundamental questions one can raise about any locally compact group}'' \cite[p.~1]{AM66}, since this property governs the structure of all its unitary representations. We refer to Subsection~\ref{subsec:prelim-type-i-groups} for precise definitions, to \cite{BH20} for a reference text on these questions, and to the introductions of \cite{Semal23} and of \cite{Mo25} (in French) for accessible presentations aimed at non-specialists.

Beyond abelian and compact groups, many classes of groups are now known to be of type~$\I$. Harish-Chandra proved in 1953 that connected semisimple Lie groups with finite center are of type~$\I$ \cite{HC53}, and Thoma proved in 1964 that a countable discrete group is of type~$\I$ only if it is virtually abelian \cite{Th64}. We refer to \cite{BH20} for further examples, and turn to totally disconnected groups, the setting of this paper. 

The first nondiscrete totally disconnected groups whose representation theory naturally attracted attention were $p$-adic groups, and more generally algebraic groups over non-Archimedean local fields. Bernstein proved in 1974 that reductive ones are of type~$\I$ \cite{Be74}. Semisimple $p$-adic groups of rank one act naturally on trees \cite{BT72}, and part of their theory soon turned out to hold for more general groups acting on trees \cite{Ti70,Ca73}, raising the question of which properties of unitary representations follow from the geometry of the action alone. Olshanskii proved that the automorphism group of a regular tree is of type~$\I$ \cite{Ol75} and classified its irreducible unitary representations into three families: spherical, special and cuspidal \cite{Ol77} (see also \cite[\S10]{Ca73}). This classification became the starting point for the study of unitary representations of closed subgroups of $\Aut(T)$, as these groups attracted growing interest, notably as a source of new examples of non-linear simple groups \cite{Ti70,BM00}.

Operator algebras play almost no role in the proofs of Bernstein and Olshanskii. In both cases, the type~$\I$ property was obtained through the following stronger, purely representation-theoretic property available for totally disconnected groups. A unitary representation $\pi:G\to\U(\H_\pi)$ of a totally disconnected locally compact group $G$ is \emph{admissible} if, for every compact open subgroup $K\leq G$,
$$
    \dim\H_\pi^K<\infty,
$$
where $\H_\pi^K$ is the space of $\pi(K)$-invariant vectors. The group $G$ is \emph{admissible} if all its irreducible unitary representations are admissible. It is \emph{uniformly admissible} if, for every compact open subgroup $K\leq G$,
$$
    \sup_\pi\dim\H_\pi^K<\infty,
$$
where $\pi$ runs over the irreducible unitary representations of $G$. These properties are related as follows:
$$
    \text{uniformly admissible}\ \Longrightarrow\ \text{admissible}\ \xLeftrightarrow{\text{Lemma~\ref{lem:ccr-admissible}}}\ \text{\textup{CCR}}\footnotemark\ \xRightarrow{\text{\cite{Ka51}}}\ \text{type~$\I$}.
$$

For\footnotetext{The definition of \textup{CCR} groups is recalled in Subsection~\ref{subsec:prelim-type-i-groups}.} every closed subgroup of $\Aut(T)$ acting transitively on the vertices and the boundary $\partial T$ of $T$, the spherical and special representations are admissible \cite{FN91}, leaving only the cuspidal ones. This led Nebbia to the \emph{type~$\I$ conjecture for groups acting on trees}, which, in its general form, predicts that every closed subgroup of $\Aut(T)$ acting two-transitively on $\partial T$, where $T$ is a thick locally finite tree, is admissible \cite[p.~311]{Ne99}, \cite[Problem~20.8.2]{CM18}, and hence of type~$\I$ \cite[Problem~20.8.3]{CM18}. We prove this conjecture with the stronger conclusion of uniform admissibility.

\begin{mainthm}
    \label{thm:main}
    Let $G$ be a closed subgroup of $\Aut(T)$ acting minimally on a locally finite tree $T$. If the action of $G$ on $\partial T$ is two-transitive, then $G$ is uniformly admissible. In particular, $G$ is of type~$\I$.
\end{mainthm}

Theorem~\ref{thm:main} was already known for groups which, in addition to two-transitivity on $\partial T$, satisfy Tits' independence property \cite{Ti70}, which underlies Olshanskii's proof of the admissibility of the cuspidal representations of $\Aut(T)$ \cite{Ol77} (see \cite{Am03,De04,Ci15}). Building on Radu's classification \cite{Ra17}, Semal extended these arguments beyond Tits' independence property, and proved Theorem~\ref{thm:main} whenever $T$ is a biregular tree whose degrees belong to an asymptotically dense subset of $\nn$ \cite[Theorems~A and~B]{Semal24}.

In contrast with these approaches, our proof of Theorem~\ref{thm:main} strongly relies on von Neumann algebras, and specifically on the group von Neumann algebra
$$
    \LL(G)=\lambda_G(G)''
$$
generated by the left regular representation $\lambda_G$ of $G$ (see Subsection~\ref{subsec:prelim-group-vna}). Let us emphasize the two main von Neumann algebraic inputs of the proof.
\begin{enumerate}[\rm (1)]
    \item \emph{$\Wstar$-amenability.} We say that $G$ is \emph{$\Wstar$-amenable} if $\LL(G)$ is amenable as a von Neumann algebra (see Subsection~\ref{subsec:prelim-amenable} for definitions). Since type~$\I$ von Neumann algebras are amenable, we use this notion as an intermediate step to prove that $\LL(G)$ is of type~$\I$. This is a genuinely operator algebraic notion, since no group-theoretic characterization of $\Wstar$-amenability is known for general locally compact groups: while amenability of $G$ implies amenability of $\LL(G)$, the converse does not hold for nondiscrete groups. 
    \item \emph{Amalgamated free products.} By Bass--Serre theory, the groups acting on trees that we are interested in can be assumed to be amalgamated free products $K_1\ast_JK_2$ of compact groups over a common open subgroup (see Subsection~\ref{subsec:prelim-bass-serre}). We then take advantage of the widely studied notion of amalgamated free product of von Neumann algebras (see Subsection~\ref{subsec:prelim-afp}), which behaves well with respect to amalgamated free products of groups, since
    $$
        \LL(K_1\ast_JK_2)\cong\LL(K_1)\ast_{\LL(J)}\LL(K_2).
    $$
    In this way, we turn our question into a purely operator algebraic one, about amalgamated free products of atomic von Neumann algebras, a setting in which many results and techniques are already available \cite{Dy93,Dy95,DR13,Re13,Ue11,Ue13}.
\end{enumerate}

These two notions were already combined by the first author and Raum in order to study the converse question, that is, to obtain geometric consequences of the type~$\I$ property for groups acting on trees. They showed that, if a closed nonamenable\footnote{Nonamenability is necessary for this direction. For instance, the affine group $\qq_p\rtimes\qq_p^\ast$ is amenable and of type~$\I$, and acts minimally on the Bruhat--Tits tree of $\operatorname{PGL}_2(\qq_p)$, but it fixes a point of the boundary.} subgroup of $\Aut(T)$ acting minimally on a thick locally finite tree $T$ is $\Wstar$-amenable, then its action satisfies a weak form of two-transitivity on $\partial T$ \cite[Theorem~C]{HR19}. This converse question had already attracted attention: in the article where he stated the conjecture solved in this paper, Nebbia proved that, if $T$ is a regular tree, then every closed unimodular vertex-transitive subgroup of $\Aut(T)$ which is admissible acts two-transitively on $\partial T$ \cite{Ne99}. Caprace and Monod asked whether the type~$\I$ property is enough \cite[Problem~20.8.3]{CM18} and, together with Kalantar, answered this question positively: a closed nonamenable subgroup of $\Aut(T)$ of type~$\I$ acting minimally on a locally finite tree $T$ acts two-transitively on $\partial T$ \cite[Corollary~D]{CKM23}. Their proof does not directly involve von Neumann algebras: it relies on a fine analysis of specific representations, the boundary representations, combined with Glimm's characterization of the type~$\I$ property in terms of weak equivalence (Theorem~\ref{thm:glimm}~\textup{(\ref*{it:glimm-weak})}). In the \hyperref[sec:amenability-boundary]{appendix}, we build on \cite[Theorem~C]{HR19} to give a von Neumann algebraic proof of their result, under the weaker assumption that the group is $\Wstar$-amenable.

Together with Theorem~\ref{thm:main}, the result of Caprace, Kalantar and Monod \cite[Corollary~D]{CKM23} yields the following theorem.

\begin{mainthm}
    \label{thm:equivalences}
    Let $G$ be a closed nonamenable subgroup of $\Aut(T)$ acting minimally on a locally finite tree $T$. The following conditions are equivalent.
    \begin{enumerate}[\rm (i)]
        \item\label{it:eq-two-transitive} The action of $G$ on $\partial T$ is two-transitive.
        \item\label{it:eq-uniform-admissible} The group $G$ is uniformly admissible.
        \item\label{it:eq-admissible} The group $G$ is admissible.
        \item\label{it:eq-ccr} The group $G$ is \textup{CCR}.
        \item\label{it:eq-type-i} The group $G$ is of type~$\I$.
        \item\label{it:eq-vn-type-i} $\LL(G)$ is of type~$\I$.
        \item\label{it:eq-vn-amenable} $\LL(G)$ is amenable \textup{(}$G$ is $\Wstar$-amenable\textup{)}.
    \end{enumerate}
\end{mainthm}

\subsection*{Outline of the proof}

In Section~\ref{sec:wstar-amenability}, we prove that if $G$ acts two-transitively on $\partial T$, then $\LL(G)$ is amenable, using a subexponential growth criterion for the amenability of operator algebras due to Kirchberg and Vaillant (Lemma~\ref{lem:ttn-subexponential-amenability}).

Section~\ref{sec:atomic-afp} is purely about von Neumann algebras. Its main result, Theorem~\ref{thm:lattice}, states that any amenable amalgamated free product of atomic von Neumann algebras is of type~$\I$, under an integrality condition on the traces of the minimal projections of the initial algebras. The proof reduces to studying the weak$^\star$ closure of an increasing union of subalgebras. Each nontrivial step is a \emph{merger} (Definition~\ref{def:merger}). Amenability then controls the type of the resulting algebra and the traces of its minimal projections.

In Section~\ref{sec:compact-amalgams}, we apply this result to amalgamated free products of compact groups. To pass from the regular representation to uniform admissibility, we use polynomial identities, in particular the Amitsur--Levitzki theorem \cite{AL50}, adapting an argument of Godement \cite[Lemma~4]{Go52}, based on work of Kaplansky \cite{Ka49} and already used by Olshanskii for tree automorphism groups \cite[Lemma~3]{Ol75}. We then deduce Theorems~\ref{thm:main} and~\ref{thm:equivalences}.

\subsection*{Origins of the project and use of AI}

The starting point of this project was our work in progress on a structure and classification theorem for amalgamated free products of atomic von Neumann algebras in the not necessarily tracial setting. We began this work on the Isle of Skye in June 2026, inspired by the work of Dykema \cite{Dy93,Dy95}, Hartglass and Nelson \cite{HN21}, Dykema and Redelmeier \cite{DR13}, and Redelmeier \cite{Re13}. By then, we were convinced that these methods could establish the type~$\I$ property of $\LL(G)$ whenever $G$ acts two-transitively on the boundary of the tree. This would provide a first step towards the type~$\I$ conjecture for groups acting on trees.

Through extensive interaction with ChatGPT-5.6 Sol and ChatGPT-6 Astra, we formulated and proved the discrete trace criterion of Theorem~\ref{thm:lattice}. ChatGPT-6 Astra supplied the argument showing that boundary two-transitivity implies $\Wstar$-amenability (Section~\ref{sec:wstar-amenability}). The argument for uniform admissibility in Section~\ref{sec:compact-amalgams} also emerged from further interaction with ChatGPT-6 Astra.

We subsequently rewrote the arguments in full and developed the definitions and conceptual framework used in this paper. We also used Claude Opus~5.5 to a lesser extent for proofreading and editing. We take full responsibility for the mathematical content of the paper.

\subsection*{Acknowledgments}

We warmly thank the organizers of the UK Operator Algebras Conference in June 2026, the International Centre for Mathematical Sciences and Sabhal M\`or Ostaig for their hospitality and for the stimulating mathematical environment on the Isle of Skye in which this project began. We also thank the Institut Henri Poincar\'e for its hospitality and the excellent working conditions provided during the 2026 thematic trimester ``Operator algebras: approximation, rigidity and dynamics''.

\tableofcontents

\section{Preliminaries}
\label{sec:preliminaries}

\subsection{Trees and their automorphisms}
\label{subsec:prelim-trees}

We follow Serre's conventions and terminology for trees \cite{Serre80}. We write $V(T)$ for the vertex set of a tree $T$ and $d$ for its combinatorial distance. We identify an \emph{oriented edge} with the ordered pair $(v,w)$ of its endpoints, and the corresponding \emph{geometric edge} with the unordered pair $\{v,w\}$. A tree is \emph{locally finite} if every vertex has finitely many neighbors.
\begin{center}
    \vspace{\medskipamount}
    \emph{All trees considered below are locally finite.}
    \vspace{\medskipamount}
\end{center}

Each tree $T$ admits a canonical bipartition of its vertex set $V(T)$ which defines two \emph{types} of vertices, with adjacent vertices of opposite types. 
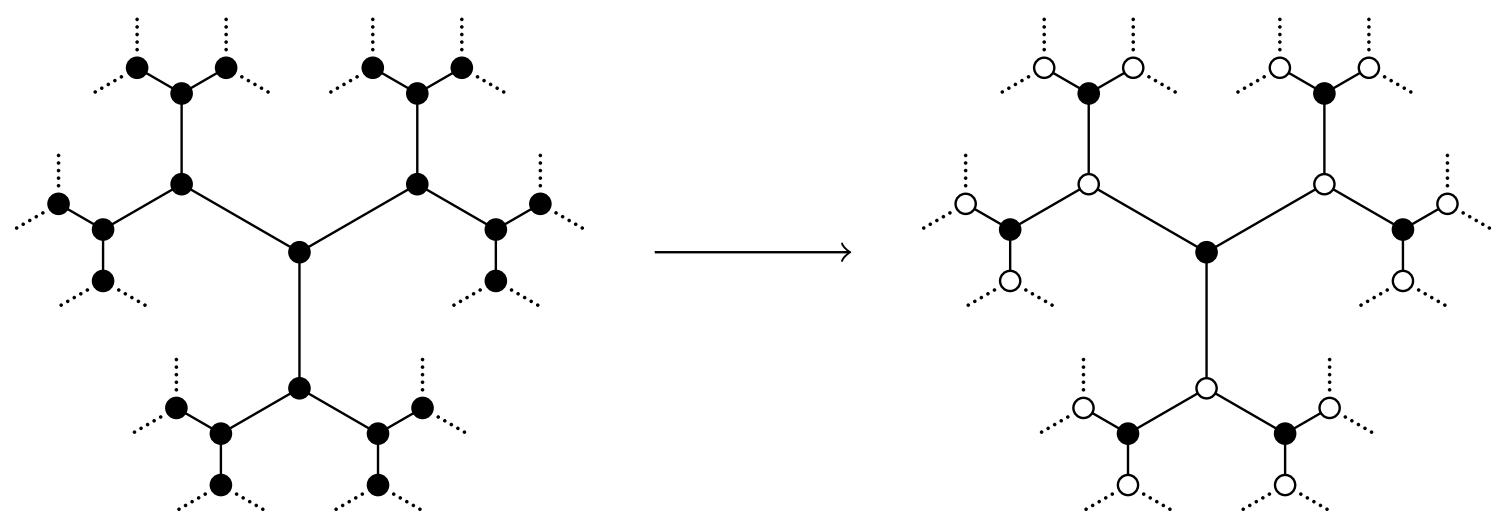
\begin{figure}[htbp]
    \centering
    \begin{tikzpicture}[vertex/.style={circle,draw,inner sep=0pt,minimum size=3.8pt},line width=0.45pt]
        \foreach \shift in {0,6} {
            \begin{scope}[xshift=\shift cm]
                \tikzset{odd vertex/.style={vertex,fill=black}}
                \ifdim\shift pt>0pt
                    \tikzset{odd vertex/.style={vertex,fill=white}}
                \fi
                \coordinate (root) at (0,0);
                \foreach \branchangle [count=\branch] in {30,150,270} {
                    \coordinate (u\branch) at (\branchangle:0.9);
                    \draw (root)--(u\branch);
                    \foreach \turn [count=\child] in {-60,60} {
                        \pgfmathsetmacro{\childangle}{\branchangle+\turn}
                        \path (u\branch) ++(\childangle:0.6) coordinate (v\branch\child);
                        \draw (u\branch)--(v\branch\child);
                        \foreach \secondturn [count=\grandchild] in {-60,60} {
                            \pgfmathsetmacro{\leafangle}{\childangle+\secondturn}
                            \path (v\branch\child) ++(\leafangle:0.34) coordinate (w\branch\child\grandchild);
                            \draw (v\branch\child)--(w\branch\child\grandchild);
                            \foreach \continuationturn in {-60,60}
                                \foreach \dotdistance in {0.12,0.17,0.22,0.27,0.32}
                                    \fill (w\branch\child\grandchild) ++({\leafangle+\continuationturn}:\dotdistance) circle[radius=0.34pt];
                            \node[odd vertex] at (w\branch\child\grandchild) {};
                        }
                        \node[vertex,fill=black] at (v\branch\child) {};
                    }
                    \node[odd vertex] at (u\branch) {};
                }
                \node[vertex,fill=black] at (root) {};
            \end{scope}
        }
        \draw[->] (2.35,0)--(3.65,0);
    \end{tikzpicture}
    \caption{A tree $T$ and its canonical bipartition of vertices}
    \label{fig:tree-bipartition}
\end{figure}

We denote by $\Aut(T)$ the group of automorphisms of a tree $T$. Equipped with the topology of pointwise convergence on $V(T)$, it is a totally disconnected second countable locally compact group, thanks to local finiteness of $T$. Any closed subgroup $G\leq \Aut(T)$ therefore inherits a locally compact group topology with the same properties. Given such a closed subgroup, we write $G_v$ for the fixator of $v\in V(T)$ and more generally $G_S$ for the fixator of any nonempty finite subset of vertices $S\subset V(T)$. These fixators are compact open subgroups of $G$ and form a neighborhood basis of the identity. Given a geometric edge $e=\{v,w\}$, we set $G_e=G_{\{v,w\}}$.

An automorphism of $T$ is \emph{type-preserving} if it preserves the type of each vertex. The type-preserving elements of a subgroup $G\leq\Aut(T)$ form a subgroup $G^+$, which is normal, open and of index at most two in $G$.

A \emph{geodesic ray} of a tree $T$ is a non-backtracking infinite sequence $(\xi_n)_{n\geq 0}$ of adjacent vertices of $T$. Two such geodesic rays are equivalent if their images differ only by finitely many vertices. The \emph{boundary} of $T$ is the set 
$$\partial T=\{\text{geodesic rays in }T\}/\sim.$$
It admits a compact topology induced by the visual metric at any vertex of $T$. An automorphism of $T$ naturally acts on geodesic rays, preserving the equivalence relation. It thus induces a homeomorphism of $\partial T$. Any closed subgroup $G\leq \Aut(T)$ therefore acts by homeomorphisms on $\partial T$: we call this action the \emph{boundary action} of $G$.

We say that $G\leq \Aut(T)$ acts \emph{minimally} on $T$ if $T$ has no proper nonempty $G$-invariant subtree. A tree is \emph{thick} if every vertex has degree at least three. We use the following reduction.

\begin{Lem}
    \label{lem:tree-thick-reduction}
    Let $G$ be a closed subgroup of $\Aut(T)$ acting minimally on $T$. Then exactly one of the following holds.
    \begin{enumerate}[\rm (i)]
        \item\label{it:tree-finite} The tree $T$ is a single vertex or a single geometric edge, and $G$ is finite.
        \item\label{it:tree-line} The tree $T$ is a bi-infinite line, and $G$ is discrete and virtually abelian.
        \item\label{it:tree-thick} The tree $T$ has a vertex of degree at least three. In this case, there is a thick tree $S$ such that $G$ is isomorphic, as a topological group, to a closed subgroup of $\Aut(S)$ acting minimally on $S$, and $\partial T$ and $\partial S$ are $G$-equivariantly homeomorphic.
    \end{enumerate}
    In particular, either $G$ is virtually abelian, or we may assume that $G$ acts minimally on a thick tree with the same boundary action.
\end{Lem}

\begin{proof}[Sketch of proof]
    This is essentially \cite[Lemma~5.2]{HR19}. If $T$ is finite, minimality forces $T$ to be a vertex or an edge, so $G$ is finite. Otherwise, minimality rules out vertices of degree one. If all vertices have degree two, $T$ is a line, and $G\leq\Aut(T)\cong D_\infty$ is discrete and virtually abelian.

    In the remaining case, $S$ is obtained from $T$ by suppressing the vertices of degree two, that is, by replacing each maximal path whose interior vertices have degree two by a single edge. Restriction to $V(S)$ identifies $\Aut(T)$ with the closed subgroup of $\Aut(S)$ preserving the lengths of these paths. This corrects \cite[Lemma~5.2]{HR19}, where the image is claimed to be all of $\Aut(S)$. Minimality and the boundary action are preserved by this construction.
\end{proof}

\subsection{Amalgamated free products of locally compact groups and Bass--Serre theory}
\label{subsec:prelim-bass-serre}

Let $G_1$ and $G_2$ be locally compact groups with a common open subgroup $J$. The \emph{amalgamated free product $G_1\ast_JG_2$} of the underlying abstract groups carries a unique group topology in which $G_1$ and $G_2$ are open subgroups with their given topologies \cite[Section~2.4]{HR19}. The neighborhoods of the identity in $J$ form a neighborhood basis of the identity. Therefore $G_1\ast_JG_2$ is a locally compact group for this topology, which is second countable when $G_1$ and $G_2$ are. If $K_1$ and $K_2$ are compact groups with a common open subgroup $J$, then $K_1\ast_JK_2$ is unimodular, since it is generated by $K_1$ and $K_2$ and the modular function is trivial on compact subgroups.

A closed subgroup $G\leq \Aut(T)$ acts \emph{without inversions} on $T$ if, for every $e\in E(T)$ and every $g\in G$,
$$ g\cdot e\neq \overline{e}.$$
Note that type-preserving automorphisms have no inversions. Let $G\leq \Aut(T)$ be a closed subgroup acting without inversions on $T$, with two orbits of vertices and one orbit of geometric edges. Then, given any edge $e=(v,w)$, we have the following isomorphism of locally compact groups \cite[Chapter~I, Section~4.1, Theorem~6]{Serre80}, \cite[Section~2.4]{HR19} 
$$ G\cong G_{v}\ast_{G_e}G_{w}.$$

We say that $G\leq \Aut(T)$ is \emph{locally transitive} if each $G_v$ acts transitively on the neighbors of $v$. When $T$ is thick, two-transitivity of $G\actson \partial T$ implies local transitivity. In fact, each $G_v$ acts transitively on every sphere about $v$ \cite[Lemma~3.1.1]{BM00}.

The following lemma is a direct consequence of Serre's theorem recalled above.

\begin{Lem}
    \label{lem:tree-amalgam}
    Let $T$ be a thick tree and $G\leq\Aut(T)$ be a closed subgroup acting locally transitively on $T$. Then, for any edge $e=(v,w)$, there is a canonical isomorphism of locally compact groups
    $$G^+\cong G_{v}\ast_{G_{(v,w)}}G_{w}.$$ 
\end{Lem}

\begin{proof}
    Automorphisms fixing a vertex are type-preserving, so $(G^+)_x=G_x$ for every $x\in V(T)$. By local transitivity, $G^+$ is transitive on the geometric edges of $T$, and it acts without inversions as it is type-preserving. The theorem of Serre recalled above \cite[Chapter~I, Section~4.1, Theorem~6]{Serre80} therefore gives $G^+=G_v\ast_{G_{(v,w)}}G_w$.
\end{proof}

\subsection{Von Neumann algebras}
\label{subsec:prelim-vna}

We refer to \cite{HR19} for the basic definitions on von Neumann algebras used in this paper, and we present our notation and conventions below. Given a Hilbert space $\H$, $\rB(\H)$ denotes the $\Cstar$-algebra of bounded operators on $\H$. 

\begin{center}
    \vspace{\medskipamount}
    \emph{In this paper, we always assume that Hilbert spaces are separable and that von Neumann algebras have separable predual.}
    \vspace{\medskipamount}
\end{center}

For a subset $S$ of some von Neumann algebra $M$, we write $\Wstar(S)$ for the von Neumann algebra generated by $S$. We write $\Z(M)$ for the center of $M$ and $\U(M)$ for its unitary group, and $M$ is a \emph{factor} if $\Z(M)=\cc1$. 
A von Neumann algebra $M$ is \emph{diffuse} if it has no minimal projections, and \emph{atomic} if every nonzero projection of $M$ dominates a minimal projection. Atomic von Neumann algebras are therefore direct sums of matrix and $\rB(\ell^2(\nn))$ blocks. 

We denote by $\csupp[M]{p}$ the \emph{central support} of a projection $p$ in a von Neumann algebra $M$, and $\csupp{p}$ when the algebra $M$ is clear from the context. A projection $p\in M$ is \emph{full} if $\csupp[M]{p}=1$. We write $\zat[M]\in\Z(M)$ for the supremum of the central supports of the minimal projections of $M$, and we put $\zdiff[M]=1-\zat[M]$. Therefore any von Neumann algebra decomposes as a direct sum of its atomic part and its diffuse part through $M=\zat[M]M\oplus \zdiff[M]M$. A projection $p\in M$ is \emph{atomic} if $p\leq\zat[M]$, equivalently if the corner $pMp$ is atomic.

A \emph{trace} on $M$ is an additive, positively homogeneous map $\tau:M_+\to[0,\infty]$ such that $\tau(x^*x)=\tau(xx^*)$ for every $x\in M$. 
\begin{center}
    \vspace{\medskipamount}
    \emph{In this paper, by a trace on a von Neumann algebra $M$ we mean a faithful normal semifinite trace on $M$. In particular, a trace is not assumed to be finite unless explicitly stated.}
    \vspace{\medskipamount}
\end{center}
If $\tau(1)<\infty$, we say that $\tau$ is a finite trace and $(M,\tau)$ is a finite tracial von Neumann algebra. A corner always carries the restriction of the ambient trace, \emph{without normalization}.

A \emph{conditional expectation} from $M$ onto a von Neumann subalgebra $N$ is a positive linear map $\rE:M\to N$ such that $\rE(b)=b$ and $\rE(bxb')=b\rE(x)b'$ for all $b,b'\in N$ and $x\in M$. An inclusion $N\subset M$ is \emph{expected} if it is unital and comes with a faithful normal conditional expectation $\rE_N^M:M\to N$. Expectations on a tower $B\subset N\subset M$ are \emph{compatible} if $\rE_B^M=\rE_B^N\circ\rE_N^M$. 

A \emph{tracial inclusion} $(N\subset M,\tau)$ is a unital inclusion of von Neumann algebras $N\subset M$ together with a trace $\tau$ on $M$ whose restriction to $N$ is semifinite. Tracial inclusions are expected: Takesaki's theorem \cite[Theorem~IX.4.2]{Ta03} gives a unique $\tau$-preserving faithful normal conditional expectation $\rE_N^M$. We will always endow tracial inclusions with this trace-preserving expectation. If $(N\subset M_1,\tau)$ and $(N\subset M_2,\tau)$ are tracial inclusions whose traces agree on $N$, we write $(M_1\supset N\subset M_2,\tau)$ and also call them tracial inclusions.

\subsection{Type \texorpdfstring{$\I$}{I} von Neumann algebras}
\label{subsec:prelim-type-i-vna}

Let $\nu$ be the Lebesgue probability measure on $[0,1]$. Then $\rL^\infty([0,1],\nu)$ is, up to isomorphism, the only abelian diffuse von Neumann algebra. It is naturally endowed with the trace obtained by integration against $\nu$. We will keep $\nu$ as a notation for this trace. Every abelian von Neumann algebra is of the form $\rL^\infty([0,1],\nu)\oplus\ell^\infty(I)\cong \rL^{\infty}(X,\mu)$ for a countable set $I$, where either summand may be zero, and $(X,\mu)$ is a measure space that may or may not have atoms.

Given $k\in\nn_{>0}$, let $\H_k$ denote the $k$-dimensional Hilbert space. The von Neumann algebra $\rB(\H_k)\cong\rM_k(\cc)$ is a factor, called the type $\I_k$ factor, which we endow with its canonical \emph{unnormalized} finite trace $\Tr_k$, so that every minimal projection has trace $1$ and $\Tr_k(1)=k$. Set $\H_{\infty}=\ell^2(\nn)$: then $\rB(\H_{\infty})$ is a factor, called the type $\I_{\infty}$ factor, which we endow with its canonical infinite trace $\Tr$ defined so that every minimal projection has trace $1$. 

The tensor product $\rB(\H_k)\ovt\rL^{\infty}(X,\mu)$ is the algebra of bounded measurable fields $x\mapsto e(x)\in\rB(\H_k)$ on $X$, and its center is then identified with $\rL^\infty(X,\mu)$.

A von Neumann algebra $M$ is \emph{of type $\I$} if 
$$
    M\cong\bigoplus_{k\in\nn_{>0}\cup\{\infty\}}\rB(\H_k)\ovt\Z_k,
$$
where $\Z_k$ are (possibly zero) abelian von Neumann algebras. 
For $t\in[0,\infty[$, a von Neumann algebra $M$ is \emph{of type $\I_{\leq t}$} if
$$
    M\cong\bigoplus_{1\leq k\leq t}\rM_k(\cc)\ovt\Z_k,
$$
where $\Z_k$ are (possibly zero) abelian von Neumann algebras.

Let $M$ be a type $\I$ von Neumann algebra, and for $k\in \nn_{>0}\cup\{\infty\}$ write $\Z_k=\rL^\infty(X_k,\mu_k)$ in the above decomposition. Viewing a projection $e\in M$ as a bounded measurable function on $X=\bigsqcup_k X_k$, taking values in $\bigsqcup_k \rB(\H_k)$, the \emph{rank function} of a projection $e\in M$ is given, for almost every $x\in X$, by
$$
    \rk_M(e)(x)=\rank e(x)\in\{0,1,2,\ldots\}\cup\{\infty\}.
$$
The rank function does not depend on the chosen decomposition: it is the center-valued dimension function of $M$ \cite[Definition~V.2.33 and Theorem~V.2.34]{Ta79}. 
We will use that the corner $eMe$ is of type $\I_{\leq t}$ exactly when $\rk_M(e)\leq t$ almost everywhere.

Finally, the following lemma characterizes type $\I_{\leq m}$ von Neumann algebras by a polynomial identity. The direct implication is the Amitsur--Levitzki theorem \cite{AL50}.

\begin{Lem}
    \label{lem:standard-identity}
    Let $M$ be a von Neumann algebra, and let $M_0\subset M$ be a weak$^\star$ dense $*$-subalgebra of $M$. For $m\in\nn_{>0}$, consider the \emph{standard polynomial} of degree $2m$ 
    $$
        S_{2m}(x_1,\ldots,x_{2m})
        =\sum_{\sigma\in\Sym_{2m}}
        \sgn(\sigma)x_{\sigma(1)}\cdots x_{\sigma(2m)}.
    $$
     Then $M$ is of type $\I_{\leq m}$ if and only if $S_{2m}$ vanishes identically on $M_0$.
\end{Lem}

\begin{proof}
    Assume that $M$ is of type $\I_{\leq m}$. Since $S_{2m}$ is multilinear and has no constant term, it suffices to show that it vanishes on $\rM_k(\cc)$ for $k\in\{1,\ldots,m\}$, which we apply fiberwise, and we may regard $\rM_k(\cc)$ as a corner of $\rM_m(\cc)$. The vanishing of $S_{2m}$ on $\rM_m(\cc)$ is the Amitsur--Levitzki theorem \cite{AL50}, and we refer to \cite{Ro76} for a short proof. In particular, $S_{2m}$ vanishes on $M_0$.

    Conversely, assume that $S_{2m}$ vanishes identically on $M_0$. By norm continuity, $S_{2m}$ vanishes on the norm closure of $M_0$. By Kaplansky's density theorem and the joint strong continuity of multiplication on bounded sets, it vanishes on $M$. Suppose that $M$ is not of type $\I_{\leq m}$. Then $M$ has a nonzero direct summand that either has no type $\I$ part or is of type~$\I_k$ with $m<k\leq\infty$. In both cases, $M$ contains matrix units $(e_{ij})_{1\leq i,j\leq m+1}$ with $e_{11}\ne0$. Evaluate $S_{2m}$ at
    $$
        e_{11},e_{12},e_{22},e_{23},\ldots,e_{mm},e_{m,m+1}.
    $$
    The displayed order is the only one giving a nonzero product, and this product is $e_{1,m+1}\ne0$. This contradiction proves the converse.
\end{proof}

\subsection{Group von Neumann algebras and Hecke algebras}
\label{subsec:prelim-group-vna}

Let $G$ be a second countable locally compact group with left Haar measure $\mu_G$ and modular function $\Delta_G$. A \emph{unitary representation} $\pi$ of $G$ is a strongly continuous homomorphism from $G$ to the unitary group of a separable Hilbert space $\H_\pi$. It integrates to a nondegenerate $*$-representation of $\rL^1(G)$ given by
$$
    \pi(f)=\int_G f(g)\pi(g)\,{\rm d}\mu_G(g),\qquad f\in\rL^1(G).
$$
Every locally compact group $G$ admits the following \emph{left regular representation} $\lambda_G$ on $\rL^2(G,\mu_G)$, given for $g\in G$, $\xi\in\rL^2(G,\mu_G)$ and almost every $x\in G$ by $$(\lambda_G(g)\xi)(x)=\xi(g^{-1}x).$$ 

Any unitary representation $(\pi,\H_\pi)$ of $G$ gives rise to the von Neumann algebra generated by its image, namely the bicommutant $\pi(G)''\subset\rB(\H_\pi)$.
In particular, the \emph{group von Neumann algebra} $\LL(G)$ of $G$ is the von Neumann algebra associated to its left regular representation: 
$$\LL(G)=\lambda_G(G)''.$$
The linear spans of $\lambda_G(G)$ and of $\lambda_G(\rC_c(G))$ are weak$^\star$ dense in $\LL(G)$, and $\lambda_G$ is injective on $\rC_c(G)$.

If $G$ is abelian with Pontryagin dual $\widehat G$, conjugating by the Fourier transform identifies $\LL(G)$ with $\rL^\infty(\widehat G)$ and $\lambda_G(f)$ with $\widehat f$ whenever $f\in \rL^1(G)$. Integration against the Haar measure $\mu_{\widehat{G}}$ then defines a trace $\tau_G$ on $\LL(G)$, which is characterized by the following \emph{Plancherel formula}: 
\begin{equation}\label{eq:plancherel}
    \tau_G(\lambda_G(f)^*\lambda_G(f))=\norm{f}_2^2, \qquad \forall f\in\rC_c(G).
\end{equation}

If $\Gamma$ is a discrete group, $\tau_{\Gamma}(x)=\ps{x\delta_e}{\delta_e}$ defines a finite trace on $\LL(\Gamma)\subset \rB(\ell^2(\Gamma))$, so that $\tau_{\Gamma}(\lambda_{\Gamma}(f))=f(e)$ for $f\in\cc[\Gamma]$. This trace still satisfies \eqref{eq:plancherel} and therefore coincides with the above trace when $\Gamma$ is abelian. For example, if $\Gamma=\zz$ we have $(\LL(\zz),\tau_{\zz})\cong (\rL^{\infty}(\mathbb{S}),\mu_{\mathbb{S}})\cong (\rL^{\infty}([0,1],\nu),\nu)$.

The two constructions above are particular cases of the general \emph{Plancherel trace} defined on $\LL(G)$ whenever $G$ is a unimodular locally compact group. It is the unique trace on $\LL(G)$ satisfying the Plancherel formula in~\eqref{eq:plancherel} \cite[Section~VII.3]{Ta03}.

Let $H\leq G$ be an open subgroup of a locally compact group $G$ and equip $H$ with the Haar measure $\mu_G|_H$. By \cite[Proposition~2.9]{HR19}, the map $\lambda_H(h)\mapsto\lambda_G(h)$ extends to a normal embedding $\LL(H)\subset\LL(G)$ with canonical faithful normal conditional expectation $\rE_{\LL(H)}^{\LL(G)}$ that satisfies 
$$
    \rE_{\LL(H)}^{\LL(G)}(\lambda_G(g))=\un_H(g)\lambda_G(g),
    \qquad g\in G.
$$
If $G$ is unimodular, then so is $H$, $\tau_G|_{\LL(H)}=\tau_H$ and $\tau_H\circ \rE_{\LL(H)}^{\LL(G)}=\tau_G$.

For a compact open subgroup $U\leq G$, put $f_U=\mu_G(U)^{-1}\un_U$ and
$$
    p_U=\lambda_G(f_U)=\frac{1}{\mu_G(U)}\int_U\lambda_G(u)\,{\rm d}\mu_G(u)\in\LL(G).
$$
It is the orthogonal projection onto the $\lambda_G(U)$-invariant vectors; we call it the \emph{averaging projection}.

The \emph{Hecke algebra}
$$
    \Hk GU=f_U\ast\rC_c(G)\ast f_U
$$
is the algebra of compactly supported $U$-biinvariant continuous functions, with convolution as product and unit $f_U$. The characteristic functions of the double cosets $UgU$ for $g\in G$ form a basis of $\Hk GU$. For $g\in G$, we have $p_U\lambda_G(g)p_U=\lambda_G(f_U\ast(L_gf_U))$, where $(L_gf_U)(x)=f_U(g^{-1}x)$ for $x\in G$. Hence $\lambda_G(\Hk GU)$ is a weak$^\star$ dense $*$-subalgebra of the \emph{Hecke von Neumann algebra} $p_U\LL(G)p_U$.

\subsection{Type \texorpdfstring{$\I$}{I} groups}
\label{subsec:prelim-type-i-groups}

Let $G$ be a second countable locally compact group. A unitary representation $(\pi,\H_{\pi})$ of $G$ is \emph{isotypical} if given any nonzero subrepresentation $\rho\leq \pi$, the multiples $\rho^{\infty}$ and $\pi^{\infty}$ are unitarily equivalent, where $\rho^{\infty}$ denotes the direct sum of countably many copies of $\rho$. Multiples of irreducible representations are obvious examples of isotypical representations. Since a representation $\pi$ is isotypical if and only if the von Neumann algebra $\pi(G)''$ it generates is a factor, they are also called \emph{factorial representations}.

The group $G$ is of \emph{type $\I$} if its isotypical representations are multiples of irreducible representations.

We record below several equivalent characterizations of the type $\I$ property, of quite different flavors. We refer to \cite[Sections~6.D and~8.F]{BH20} for precise definitions, statements and proofs, the difficult implications being mostly due to Glimm \cite{Gl61}. Condition~\textup{(\ref*{it:glimm-vna})} justifies the ``type $\I$'' terminology. The equivalence with condition~\textup{(\ref*{it:glimm-disintegration})} says that $G$ is of type $\I$ if and only if it has a ``tame'' representation theory, in which unitary representations can be recovered from their irreducible components, as in the case of abelian and compact groups.

\begin{The}[{\cite{Gl61}, \cite[Theorem~8.F.3]{BH20}}]
    \label{thm:glimm}
    Let $G$ be a second countable locally compact group. The following conditions are equivalent.
    \begin{enumerate}[\rm (i)]
        \item $G$ is a type $\I$ group.
        \item\label{it:glimm-vna} For every representation $\pi$ of $G$, $\pi(G)''$ is a type $\I$ von Neumann algebra.
        \item\label{it:glimm-disintegration} Every representation of $G$ disintegrates as a direct integral of irreducible unitary representations in an essentially unique way.
        \item\label{it:glimm-weak} Any two weakly equivalent irreducible representations of $G$ are unitarily equivalent.
        \item\label{it:glimm-dual} The topological and Borel structures of the unitary dual $\widehat{G}$ are well behaved.
    \end{enumerate}
\end{The}

Compact and abelian groups are the archetypes of type $\I$ groups. In contrast, discrete groups are archetypes of non-type $\I$ groups, as a discrete group is of type $\I$ only if it is virtually abelian \cite{Th64}. 

Apart from these, establishing the type $\I$ property is usually difficult, and often amounts to classifying all irreducible unitary representations of a group. Dixmier proved that connected nilpotent Lie groups and connected real algebraic groups are of type $\I$, while Harish-Chandra and Bernstein proved that reductive algebraic groups over Archimedean and non-Archimedean local fields, respectively, are of type $\I$ \cite[Theorem~6.E.19]{BH20}. Olshanskii proved that the full automorphism group of a regular tree is of type $\I$ \cite{Ol75,Ol77}.

A locally compact group $G$ is \emph{\textup{CCR}} if, for every irreducible unitary representation $\pi$ of $G$, $\pi(\rC_c(G))$ consists of compact operators. By a theorem of Kaplansky, \textup{CCR} groups are of type $\I$ \cite{Ka51}. This condition is especially useful to prove the type $\I$ property when the group under consideration is totally disconnected, as the following characterization shows. 

Assume that $G$ is totally disconnected, and, for every representation $(\pi,\H_{\pi})$ and compact open subgroup $U$, let $\H_{\pi}^{U}$ be the space of $\pi(U)$-invariant vectors of $\H_\pi$. Note that $\H_{\pi}^{U}=\pi(f_U)\H_\pi$. Recall that $G$ is \emph{admissible} if $\dim\H_\pi^U<\infty$ for every irreducible representation $\pi$ and every $U$, and \emph{uniformly admissible} if moreover, for every $U$, these dimensions are bounded independently of $\pi$. 

\begin{Lem}
    \label{lem:ccr-admissible}
    A totally disconnected locally compact group is \textup{CCR} if and only if it is admissible. In particular, uniformly admissible groups are \textup{CCR}.
\end{Lem}

\begin{proof}
    Let $\pi$ be an irreducible unitary representation of a \textup{CCR} group $G$. Then, for every compact open subgroup $U\leq G$, the projection $\pi(f_U)$ onto $\H_\pi^U$ is a compact operator, hence of finite rank, so that $\H_\pi^U$ is finite dimensional. Conversely, assume that every $\H_\pi^U$ is finite dimensional, and let $f\in\rC_c(G)$. For every compact open subgroup $U\leq G$, the operator $\pi(f_U\ast f\ast f_U)=\pi(f_U)\pi(f)\pi(f_U)$ has finite rank. As $U$ runs through the compact open subgroups of $G$, ordered by reverse inclusion, $f_U\ast f\ast f_U\to f$ in $\rL^1(G)$, so that $\pi(f_U\ast f\ast f_U)\to\pi(f)$ in norm. Hence $\pi(f)$ is compact.
\end{proof}

\subsection{Amenable von Neumann algebras, $\Wstar$-amenable groups}
\label{subsec:prelim-amenable}

A von Neumann algebra $M$ is \emph{amenable} if, for some faithful normal unital representation $M\subset\rB(\H)$, there exists a (not necessarily normal) conditional expectation $\rE:\rB(\H)\to M$. 

We first record some classical facts about amenable von Neumann algebras. 

\begin{fact}[{\cite[Sections~III.3 and IV.2]{Bl06}}]
    \label{fact:amenability-properties}
    Let $M$ be a von Neumann algebra.
    \begin{enumerate}[\rm (i)]
        \item\label{it:amenability-expected-corners} Assume that $M$ is amenable. Then for any projection $p\in M$, $pMp$ is amenable. If $\rE:M\to N$ is a conditional expectation, then $N$ is amenable. 
        \item\label{it:amenability-increasing-corners} 
        If $(p_j)_j$ is an increasing net of projections converging weakly$^\star$ to $1$ and each $p_jMp_j$ is amenable, then $M$ is amenable.
        \item\label{it:amenability-type-i} Type $\I$ von Neumann algebras are amenable.
    \end{enumerate}
\end{fact}

We say that a locally compact group $G$ is \emph{$\Wstar$-amenable} if $\LL(G)$ is amen\-able: this is the central notion of this paper. By the results below, this property coincides with amenability for discrete groups, but is strictly weaker in general. To date, it admits no purely group-theoretic characterization.
\begin{fact} Let $G$ be a locally compact group.
    \begin{enumerate}[\rm (i)]
        \item If $G$ is amenable as a group, then $\LL(G)$ is an amenable von Neumann algebra.
        \item If $G$ is discrete and $\LL(G)$ is amenable, then $G$ is amenable. 
        \item If $G$ is connected and second countable, then $\LL(G)$ is amenable.
    \end{enumerate}
\end{fact}
\begin{proof}
\begin{enumerate}[\rm (i)]
    \item Let $m$ be a left-invariant mean on $\rL^\infty(G)$ and let $\rho$ denote the right regular representation of $G$. Then
    $$
        \ps{\rE(T)\xi}{\eta}=m\big(g\mapsto\ps{\rho(g)T\rho(g)^*\xi}{\eta}\big)
    $$
    defines a conditional expectation $\rE:\rB(\rL^2(G))\to\rho(G)'=\LL(G)$.
    \item Assume that $G$ is discrete, and let $\rE:\rB(\ell^2(G))\to\LL(G)$ be a conditional expectation onto $\LL(G)$. Then $f\mapsto\tau_G(\rE(M_f))$ is a left-invariant mean on $\ell^\infty(G)$, $M_f$ being the multiplication operator by $f$.
    \item This is a corollary of the deep classification of amenable factors due to Connes \cite[Corollary~6.9]{Co76}.\qedhere
\end{enumerate}
\end{proof}

Many results on amenable groups have operator algebraic counterparts (see Proposition~\ref{prop:scalar-amenability} below). For instance, the following result, due to Kirchberg and Vaillant (in a stronger $\Cstar$-algebraic form), is the von Neumann algebraic analogue of the fact that groups of subexponential growth are amenable.

\begin{Lem}[{\cite[Corollary~2.2]{KV92}}]
    \label{lem:ttn-subexponential-amenability}
    Let $M$ be a von Neumann algebra. Suppose that $M$ has finite-dimensional linear subspaces $(\A_n)_{n\geq0}$ satisfying, for all $n,m\in\nn$,
    $$
        \begin{gathered}
            1\in\A_0,\qquad\A_n^*=\A_n,\qquad\A_n\subset\A_{n+1},\qquad\A_n\A_m\subset\A_{n+m},\\
            M=\Wstar\Bigl(\bigcup_{n\geq0}\A_n\Bigr),\qquad\text{and}\qquad\limsup_{n\to\infty}\frac{\log\dim\A_n}{n}=0.
        \end{gathered}
    $$
    Then $M$ is amenable.
\end{Lem}

\begin{proof}
    The norm closure $\A$ of $\bigcup_{n\geq0}\A_n$ is a nuclear $\Cstar$-algebra \cite[Corollary~2.2]{KV92}. Since $\A$ is weak$^\star$ dense in $M$, the algebra $M$ is amenable \cite[Theorem~IV.2.2.13]{Bl06}.
\end{proof}

\subsection{Amalgamated free products of von Neumann algebras}
\label{subsec:prelim-afp}

Let $A_1\supset B\subset A_2$ be expected inclusions of von Neumann algebras. For $i\in\{1,2\}$, an element $x\in A_i$ is \emph{centered} if $\rE_B^{A_i}(x)=0$, and a \emph{reduced word} is a product $x_1\cdots x_n$ of centered elements $x_j\in A_{i_j}$, where $n\in\nn_{>0}$ and $i_1,\ldots,i_n\in\{1,2\}$, with $i_j\neq i_{j+1}$ for $1\leq j<n$. 

The \emph{amalgamated free product} $M=A_1\ast_B A_2$ is the von Neumann algebra, equipped with a faithful normal conditional expectation $\rE_B^M\colon M\to B$, such that:
\begin{itemize}
    \item $M$ is generated by copies of $A_1$ and $A_2$ containing a common copy of $B$;
    \item $\rE_B^M$ extends $\rE_B^{A_1}$ and $\rE_B^{A_2}$;
    \item every reduced word $x_1\cdots x_n$ satisfies $\rE_B^M(x_1\cdots x_n)=0$.
\end{itemize}
It exists, and it is unique up to an isomorphism that is the identity on $A_1$ and $A_2$ and intertwines the expectations \cite{Vo85}, see also \cite{Ue99}. In particular, $A_1\cap A_2=B$ inside $M$.
If $(A_1\supset B\subset A_2,\tau)$ is a tracial inclusion, then $\tau\circ\rE_B^M$ is a trace on $M=A_1\ast_BA_2$ extending $\tau$, and we always equip $M$ with it. Finally, when $B=\cc1$, we recover the classical notion of free product of von Neumann algebras, which we denote by $A_1\ast A_2$. In this case, the expectations are faithful normal states $\phi_i$ on $A_i$.

Amalgamated free products of von Neumann algebras are compatible with those of groups: if $K_1$ and $K_2$ are second countable locally compact groups with a common open subgroup $J$, then by \cite[Proposition~2.22]{HR19}, with the canonical expectations of the open subgroups,
\begin{equation}\label{eq:compact-amalgam-afp}
    \LL(K_1\ast_{J}K_2)\cong\LL(K_1)\ast_{\LL(J)}\LL(K_2).
\end{equation}

As for groups, where $\zz/2\ast\zz/2$ is the only amenable free product of two nontrivial groups, amenability of a free product of von Neumann algebras is very restrictive, as shown by the following result of Ueda.

\begin{Prop}[{\cite[Remarks~4.2(5)]{Ue11}}]
    \label{prop:scalar-amenability}
    Let $(A_1,\tau_1)$ and $(A_2,\tau_2)$ be von Neumann algebras with finite traces satisfying $\tau_1(1_{A_1})=\tau_2(1_{A_2})$, and assume that $A_1,A_2\neq\cc1$. Then $A_1\ast A_2$ is amenable if and only if $A_1\cong A_2\cong\cc^2$.
\end{Prop}

In this exceptional case, the free product can be computed explicitly, and it is of type $\I$ with matrix sizes at most~$2$, as shown by the following result of Dykema.
Here and below, labels below summands specify their traces, and summands with zero trace are omitted.

\begin{samepage}
\begin{Prop}[{\cite[Theorem~1.1]{Dy93}}]
    \label{prop:two-projections}
    Let $\alpha,\beta,\gamma\in\rr$ satisfy $\gamma>\alpha\geq\max\{\beta,\gamma-\beta\}>0$. Then
    $$
        \left(\underset{\alpha}{\cc}\oplus\underset{\gamma-\alpha}{\cc}\right)
        \ast\left(\underset{\beta}{\cc}\oplus\underset{\gamma-\beta}{\cc}\right)
        \cong\underset{\alpha+\beta-\gamma}{\cc}
        \oplus\underbrace{\bigl(\rM_2(\cc)\ovt\rL^\infty([0,1],\nu)\bigr)}_{(\gamma-\alpha)(\Tr_2\pt\nu)}
        \oplus\underset{\alpha-\beta}{\cc}.
    $$
    Moreover, on the diffuse part of the free product, every minimal projection of either free factor has rank one almost everywhere.
\end{Prop}
\end{samepage}

In particular, an amenable free product of two nontrivial von Neumann algebras is automatically of type $\I$, and the minimal projections of both free factors have rank at most one in it. Section~\ref{sec:atomic-afp} builds on these two results to establish an amalgamated version of this phenomenon, Theorem~\ref{thm:lattice}, in which the ranks of projections are controlled by their traces.

\section{From boundary two-transitivity to \texorpdfstring{$\Wstar$}{W*}-amenability}
\label{sec:wstar-amenability}

In this section, we prove the following result.

\begin{Prop}
    \label{prop:ttn-boundary-amenability}
    Let $T$ be a thick locally finite tree and let $G$ be a closed subgroup of $\Aut(T)$. If $G$ acts two-transitively on $\partial T$, then $G$ is $\Wstar$-amenable.
\end{Prop}

The main step relates spherical transitivity to the growth of Hecke algebras. Under the hypothesis below, the number of double cosets of bounded displacement grows at most linearly, allowing us to apply Lemma~\ref{lem:ttn-subexponential-amenability} to the corresponding finite corners.

\begin{Prop}
    \label{prop:ttn-spherical-amenability}
    Let $G$ be a closed subgroup of $\Aut(T)$, where $T$ is a locally finite tree. Fix $v\in V(T)$ and put $K=G_v$. Suppose that, for every $n\in\nn$, the group $K$ acts transitively on
    $$
        Gv\cap\{w\in V(T)\mid d(v,w)=n\}
    $$
    whenever this set is nonempty. Then $G$ is $\Wstar$-amenable.
\end{Prop}

\begin{proof}
    The subgroup $K$ is compact and open. Put $\ell(g)=d(v,gv)$ for $g\in G$. The function $\ell$ is continuous, integer valued, constant on the double cosets of $K$, and satisfies, for all $g,h\in G$,
    $$
        \ell(g^{-1})=\ell(g),\qquad\ell(gh)\leq\ell(g)+\ell(h).
    $$
    The map $gK\mapsto gv$ identifies $G/K$ with $Gv$. If $\ell(g)=\ell(h)$, the transitivity assumption gives $hv=kgv$ for some $k\in K$, so $(kg)^{-1}h\in K$ and $h\in KgK$, showing that each nonempty length level is a single double coset of $K$. Consequently, the assumption on spheres implies that, for every $n\in\nn$, the set $\{g\in G\mid\ell(g)\leq n\}$ is a union of at most $n+1$ double cosets of $K$.

    Let $U\leq K$ be a compact open subgroup and put $r=[K:U]\in\nn_{>0}$. Recall the projection $p_U=\lambda_G(f_U)$, and put
    $$
        N_U=p_U\LL(G)p_U.
    $$
    For $n\in\nn$, let $\Hk GU_n$ be the subspace of $\Hk GU$ consisting of functions supported in $\{g\in G\mid\ell(g)\leq n\}$. The characteristic functions of the double cosets of $U$ contained in this set form a basis of $\Hk GU_n$. These double cosets are compact and open. Choose $a_1,\ldots,a_r,b_1,\ldots,b_r\in K$ such that $K=\bigcup_{i=1}^rUa_i=\bigcup_{j=1}^rb_jU$. For $g\in G$, we have $KgK=\bigcup_{i,j=1}^rUa_igb_jU$, where some of these double cosets may coincide. Thus each double coset of $K$ is a union of at most $r^2$ double cosets of $U$, and
    $$
        \dim\Hk GU_n\leq r^2(n+1).
    $$
    For the involution and convolution of $\rC_c(G)$, we have, for all $n,m\in\nn$,
    $$
    \begin{gathered}
        \Hk GU_n^*=\Hk GU_n,\\
        \Hk GU_n\ast\Hk GU_m\subset\Hk GU_{n+m}.
    \end{gathered}
    $$
    The first assertion follows from $\ell(g^{-1})=\ell(g)$, since $\Delta_G=1$ on compact subgroups. For the second, let $f_1\in\Hk GU_n$, $f_2\in\Hk GU_m$ and $g\in G$. If $(f_1\ast f_2)(g)\ne0$, then $g=ab$ with $a\in\supp(f_1)$ and $b\in\supp(f_2)$, so that $\ell(g)\leq\ell(a)+\ell(b)\leq n+m$.

    Put $\A_{U,n}=\lambda_G(\Hk GU_n)$ for every $n\in\nn$. These subspaces form an increasing self-adjoint multiplicative filtration of $N_U$, with $p_U\in\A_{U,0}$ and $\dim\A_{U,n}\leq r^2(n+1)$. Their union is weak$^\star$ dense in $N_U$ by the end of Subsection~\ref{subsec:prelim-group-vna}. Lemma~\ref{lem:ttn-subexponential-amenability} shows that $N_U$ is amenable.

    Let $(U_j)_{j\geq1}$ be the fixators in $G$ of the balls of radius $j$ about $v$. They form a decreasing neighborhood basis of the identity consisting of compact open subgroups of $K$, and the projections $p_j=p_{U_j}$ satisfy $p_j\leq p_{j+1}$. Strong continuity of $\lambda_G$ gives, for every $\xi\in\rL^2(G)$,
    $$
        \norm{p_j\xi-\xi}
        \leq\sup_{u\in U_j}\norm{\lambda_G(u)\xi-\xi}\longrightarrow0.
    $$
    Thus $p_j\to1$ strongly, hence in the weak$^\star$ topology.

    Since each corner $p_j\LL(G)p_j$ is amenable, so is $\LL(G)$ by Fact~\ref{fact:amenability-properties}\textup{(\ref*{it:amenability-increasing-corners})}. Thus $G$ is $\Wstar$-amenable.
\end{proof}

To deduce the boundary statement, it remains to obtain the spherical transitivity assumption from the action on $\partial T$.

\begin{proof}[Proof of Proposition~\ref{prop:ttn-boundary-amenability}]
    By \cite[Lemma~3.1.1]{BM00}, each vertex stabilizer $G_v$ acts transitively on every sphere about $v$. Proposition~\ref{prop:ttn-spherical-amenability} therefore applies and shows that $G$ is $\Wstar$-amenable.
\end{proof}

\section{From amenability to type \texorpdfstring{$\I$}{I} for AFP von Neumann algebras}
\label{sec:atomic-afp}

In this section, we prove that an amenable amalgamated free product of atomic von Neumann algebras is of type~$\I$, provided the traces of the minimal projections lie in $\delta\nn$ for some $\delta>0$.

\begin{The}
    \label{thm:lattice}
    Let $(A_1\supset B\subset A_2,\tau)$ be tracial inclusions of atomic von Neumann algebras. Assume that there exists some $\delta\in\rr_{>0}$ such that, for each $i\in\{1,2\}$ and every minimal projection $p\in A_i$,
    \begin{equation}\label{eq:lattice-hypothesis}
        \tau(p)\in\delta\nn.
    \end{equation}
    If $M=A_1\ast_B A_2$ is amenable, then $M$ is of type~$\I$. Moreover, $pMp$ is of type $\I_{\leq\tau(p)/\delta}$ for every finite-trace projection $p\in A_1\cup A_2$.
\end{The}

When $A_1$ and $A_2$ are finite dimensional, condition~\eqref{eq:lattice-hypothesis} is not needed: if $M$ is amenable, then $M$ is of type~$\I$ by \cite[Theorem~5.1]{Dy95}. In general, however, it cannot be dropped: when the traces of the minimal projections of $A_1$ and $A_2$ come arbitrarily close to zero, $M$ can be the hyperfinite type $\II_1$ factor \cite[Example~5.5]{Dy95}.

Note that every minimal projection $p\in B$ also satisfies $\tau(p)\in\delta\nn$. Indeed, semifiniteness on $B$ gives $\tau(p)<\infty$, and, viewed in $A_1$, the projection $p$ is a finite orthogonal sum of minimal projections, each of trace in $\delta\nn$.

\subsection{Classical results on AFP von Neumann algebras}

{\emergencystretch=1em
We recall the following well-known identities for amalgamated free products of von Neumann algebras.\par}

\begin{Lem}
    \label{lem:afp-identities}
    Let $A_1\supset B\subset A_2$ be expected inclusions of von Neumann algebras.
    \begin{enumerate}[\rm (i)]
        \item\label{it:afp-full-corner} If $p\in B$ is a projection with $\csupp[B]{p}=1$, then there is an isomorphism, compatible with the compressed expectations, such that
        $$
            p(A_1\ast_B A_2)p\cong pA_1p\ast_{pBp}pA_2p.
        $$
        \item\label{it:afp-compression} If $r\in\Z(B)$ is a nonzero projection such that $A_2=rA_2r\oplus B(1-r)$, then there is an isomorphism, compatible with the compressed expectations, such that
        $$
            r(A_1\ast_B A_2)r\cong rA_1r\ast_{rB}rA_2r.
        $$
        \item\label{it:afp-refinement} If $B\subset\tilde{B}\subset A_1$ is an intermediate von Neumann algebra with compatible expectations, then
        $$
            A_1\ast_B A_2\cong A_1\ast_{\tilde{B}}(\tilde{B}\ast_B A_2).
        $$
    \end{enumerate}
\end{Lem}

\begin{proof}
    \textup{(\ref*{it:afp-full-corner})} This is \cite[Lemma~4.6]{Ue13}.

    \textup{(\ref*{it:afp-compression})} This follows from \cite[Lemma~4.3]{Dy95}, whose proof also applies in the non-tracial setting. See also \cite{Ue16}.

    \textup{(\ref*{it:afp-refinement})} The reduced $\Cstar$-algebra version is \cite[Corollary~3.9]{NSS02}. Choosing a faithful normal state on $B$ and composing it with the expectations, the state-preserving isomorphism extends to the von Neumann algebra closures in the corresponding GNS representations.
\end{proof}

\begin{fact}
    \label{fact:afp-subalgebras}
    Let $(\widetilde A_1\supset B\subset\widetilde A_2,\tau)$ be a tracial inclusion, and let $B\subset A_i\subset\widetilde A_i$ be von Neumann subalgebras for $i\in\{1,2\}$. Then the von Neumann subalgebra of $\widetilde A_1\ast_B\widetilde A_2$ generated by $A_1$ and $A_2$ is identified with $A_1\ast_B A_2$, and this inclusion is trace-preserving, hence expected.
\end{fact}

Indeed, the trace-preserving expectation of $A_i$ onto $B$ is the restriction of that of $\widetilde A_i$, so reduced words in $A_1\ast_BA_2$ are reduced words in $\widetilde A_1\ast_B\widetilde A_2$. The conclusion now follows from the uniqueness of the amalgamated free product.

\subsection{Mergers and type \texorpdfstring{$\I$}{I} stability}

In this subsection, we introduce the notion of \emph{mergers}, which morally describe the addition of a partial isometry intertwining two orthogonal projections of equal trace of a given von Neumann algebra. We will use this operation to approximate a given AFP $A_1\ast_D A_2$, where, for each $i\in\{1,2\}$, the atomic von Neumann algebra $A_i$ is generated by an atomic maximal abelian self-adjoint subalgebra (masa) $D_i\supset D$ and matrix units in $A_i$. Adjoining these matrix units one at a time realizes each $A_i$ as the weak$^\star$ closure of an increasing union of successive mergers of $D_i$ (Example~\ref{ex:block-merger}). By Lemma~\ref{lem:outer-fusion}, $A_1\ast_D A_2$ is then the weak$^\star$ closure of an increasing union of successive mergers of $D_1\ast_D D_2$.\footnote{The idea behind mergers is not new. After reducing to an abelian amalgam \cite[Lemma~5.2]{Dy95}, Dykema adjoins matrix units $v$, called \emph{connectors}, one at a time \cite[pp.~1583--1585]{Dy95}, which amounts to taking the mergers $M_{vv^*,v^*v}$ of Definition~\ref{def:merger}. Merging two matrix blocks also appears as the second type of simple step in \cite[Definition~3.4(2)]{DR13}. 
} 

\begin{Def}
    \label{def:merger}
    Let $(M,\tau)$ be a von Neumann algebra with a trace. Let $p,q\in M$ be orthogonal projections with $0<\tau(p)=\tau(q)<\infty$ and put $s=1-(p+q)$. Embed $\cc p\oplus\cc q\oplus sMs$ diagonally in $\rM_2(\cc)\oplus sMs$.

    The \emph{merger of $p$ and $q$ in $M$} is the von Neumann algebra
    \begin{equation}\label{eq:free-merger}
        M_{p,q}=M\ast_{\cc p\oplus\cc q\oplus sMs}\bigl(\rM_2(\cc)\oplus sMs\bigr),
    \end{equation}
    where the expectations are the trace-preserving ones.
\end{Def}

The von Neumann algebra $M_{p,q}$ is generated by $M$ and the matrix unit $v\in\rM_2(\cc)$ with $vv^*=p$ and $v^*v=q$, which we call the \emph{adjoined matrix unit}. By Subsection~\ref{subsec:prelim-afp}, $M_{p,q}$ carries a trace extending $\tau$, which we still denote by $\tau$.

\begin{Ex}
    \label{ex:block-merger}
    Let $k,l\in\nn_{>0}$ and consider a von Neumann algebra with a trace
    $$
        (A,\tau)=(A_0,\tau_{\mid A_0})\oplus(\rM_k(\cc),\Tr_k)\oplus(\rM_l(\cc),\Tr_l),
    $$ and let $p$ and $q$ be minimal projections respectively in the $\rM_k(\cc)$ and $\rM_l(\cc)$ summands.

    There is a trace-preserving isomorphism, which is the identity on $A$,
    $$
        (A_{p,q},\tau)\cong(A_0,\tau_{\mid A_0})\oplus(\rM_{k+l}(\cc),\Tr_{k+l}).
    $$
\end{Ex}
\begin{proof}
    Let $v\in A_{p,q}$ be the adjoined matrix unit, and choose matrix units $(f_{ij})$ and $(g_{ab})$ of the two matrix summands with $f_{11}=p$ and $g_{11}=q$. Since $v=pvq$, the projection $1_{A_0}$ is central in $A_{p,q}$, and the elements $f_{ij}$, $g_{ab}$, $f_{i1}vg_{1a}$ and $(f_{i1}vg_{1a})^*$ form a system of matrix units of size $k+l$ with unit $1-1_{A_0}$. Together with $A_0$, they generate $A_{p,q}$, and the trace of $A_{p,q}$ takes the value $1$ on $p$.
\end{proof}

We now show that mergers are compatible with amalgamated free products: merging two projections of one initial algebra amounts to merging them in the free product.

\begin{Lem}
    \label{lem:outer-fusion}
    Let $(A_1\supset D\subset A_2,\tau)$ be tracial inclusions of von Neumann algebras. Let $p,q\in A_1$ be orthogonal projections with $0<\tau(p)=\tau(q)<\infty$, put $s=1-(p+q)$, and assume that $D\subset\cc p\oplus\cc q\oplus sA_1s$. Then
    \[
        (A_1)_{p,q}\ast_{D}A_2\cong \bigl(A_1\ast_{D}A_2\bigr)_{p,q}.
    \]
    The isomorphism is trace-preserving and is the identity on $A_1\ast_{D}A_2$.
\end{Lem}

\begin{proof}
    All expectations are trace-preserving, hence compatible. Put $B=\cc p\oplus\cc q\oplus sA_1s$ and $\widetilde{B}=\cc p\oplus\cc q\oplus s(A_1\ast_DA_2)s$. Then, using Lemma~\ref{lem:afp-identities}\textup{(\ref*{it:afp-refinement})} with the inclusions $D\subset B\subset A_1$ and $B\subset\widetilde{B}\subset A_1\ast_DA_2$ together with associativity, we obtain
    \begin{align*}
        (A_1)_{p,q}\ast_DA_2&
        =\bigl(A_1\ast_{B}(\rM_2(\cc)\oplus sA_1s)\bigr)\ast_DA_2\\
        &\cong\bigl(A_1\ast_{B}(\rM_2(\cc)\oplus sA_1s)\bigr)\ast_{B}\bigl(B\ast_DA_2\bigr)\\
        &\cong\Bigl(A_1\ast_{B}\bigl(B\ast_DA_2\bigr)\Bigr)\ast_{B}(\rM_2(\cc)\oplus sA_1s)\\
        &\cong(A_1\ast_DA_2)\ast_{B}(\rM_2(\cc)\oplus sA_1s)\\
        &\cong(A_1\ast_DA_2)\ast_{\widetilde{B}}\Bigl(\widetilde{B}\ast_{B}(\rM_2(\cc)\oplus sA_1s)\Bigr)\\
        &\cong(A_1\ast_DA_2)\ast_{\widetilde{B}}\bigl(\rM_2(\cc)\oplus s(A_1\ast_DA_2)s\bigr)=(A_1\ast_DA_2)_{p,q}.
        \end{align*}
 The last isomorphism is obtained as follows. Since $s$ is central in $B$, $\widetilde{B}$ and $\rM_2(\cc)\oplus sA_1s$, it is central in $\widetilde{B}\ast_{B}(\rM_2(\cc)\oplus sA_1s)$. As $\widetilde{B}=s\widetilde{B}s\oplus B(1-s)$ and $\rM_2(\cc)\oplus sA_1s=(1-s)\bigl(\rM_2(\cc)\oplus sA_1s\bigr)(1-s)\oplus Bs$, Lemma~\ref{lem:afp-identities}\textup{(\ref*{it:afp-compression})}, applied to $1-s$ and $s$, gives
    \begin{align*}
        (1-s)\bigl(\widetilde{B}\ast_{B}(\rM_2(\cc)\oplus sA_1s)\bigr)&\cong(\cc p\oplus\cc q)\ast_{\cc p\oplus\cc q}\rM_2(\cc)=\rM_2(\cc),\\
        s\bigl(\widetilde{B}\ast_{B}(\rM_2(\cc)\oplus sA_1s)\bigr)&\cong s(A_1\ast_DA_2)s\ast_{sA_1s}sA_1s=s(A_1\ast_DA_2)s.
    \end{align*}
    Hence $\widetilde{B}\ast_{B}(\rM_2(\cc)\oplus sA_1s)\cong\rM_2(\cc)\oplus s(A_1\ast_DA_2)s$.
\end{proof}

When $p+q=1$, a merger is determined by its corner $pM_{p,q}p$. Following \cite[Lemma~4.3]{Dy95} and \cite[Lemmas~2.4 and~2.7]{Re13}, we compute this corner in two cases.

\begin{Lem}
    \label{lem:equivalent-endpoints}
    Let $(M,\tau)$ be a von Neumann algebra with a finite trace and let $p,q\in M$ be projections with $p+q=1$ and $0<\tau(p)=\tau(q)$. Then $M_{p,q}\cong\rM_2(\cc)\ovt pM_{p,q}p.$ Moreover:
    \begin{enumerate}[\rm (i)]
        \item\label{it:separated-corners} If $pMq=0$, then
        $$
            pM_{p,q}p\cong pMp\ast qMq.
        $$
        \item\label{it:equivalent-corners} If $p$ and $q$ are Murray--von Neumann equivalent in $M$, then
        $$
            pM_{p,q}p\cong pMp\ast\LL(\zz),
        $$
        where $\LL(\zz)$ is equipped with the trace $\tau(p)\tau_\zz$.
    \end{enumerate}
\end{Lem}

\begin{proof}
    Put $D=\cc p\oplus\cc q$, and let $v\in M_{p,q}$ be the adjoined matrix unit. Since $p,v,v^*,q$ are $2\times2$ matrix units in $M_{p,q}$ with $p+q=1$, the map
    $$
        M_{p,q}\longrightarrow\rM_2(\cc)\ovt pM_{p,q}p,\qquad
        x\longmapsto\begin{pmatrix}pxp&pxv^*\\vxp&vxv^*\end{pmatrix}
    $$
    is an isomorphism.

    \textup{(\ref*{it:separated-corners})} Since $pMq=0$, we have $M=pMp\oplus qMq$. The algebra $M_{p,q}$ is generated by $M$ and $v=pvq$. In a nonzero product $p\,a_1\cdots a_k\,p$ with $a_1,\ldots,a_k\in pMp\cup qMq\cup\{v,v^*\}$, the letters from $qMq$ occur in blocks $v\,b_1\cdots b_m\,v^*$ with $b_1,\ldots,b_m\in qMq$. Hence $pM_{p,q}p=\Wstar(pMp,\,v(qMq)v^*)$. Moreover, $y\mapsto vyv^*$ is a trace-preserving isomorphism from $qMq$ onto $v(qMq)v^*$. Let $a\in pMp$ and $b\in qMq$ with $\tau(a)=\tau(b)=0$. Then $\rE_D^M(a)=\rE_D^M(b)=0$, while $\rE_D^{\rM_2(\cc)}(v)=\rE_D^{\rM_2(\cc)}(v^*)=0$. Hence, once each $vbv^*$ is expanded, every alternating product of such elements $a$ and $vbv^*$ is an alternating product of elements of $\ker\rE_D^M$ and $\ker\rE_D^{\rM_2(\cc)}$, that is, a reduced word in $M_{p,q}$. By freeness with amalgamation over $D$, its image under $\rE_D^{M_{p,q}}$ is zero, so it has trace zero. Thus $pMp$ and $v(qMq)v^*$ are free in $pM_{p,q}p$, which proves~\textup{(\ref*{it:separated-corners})}.

    \textup{(\ref*{it:equivalent-corners})} Choose $u\in M$ with $uu^*=p$ and $u^*u=q$, so that $u=puq$, and put $h=vu^*\in\U(pM_{p,q}p)$. Since $M$ is generated by $pMp$ and $u$, and $v=hu$, the algebra $M_{p,q}$ is generated by $pMp$, $h$ and $u$. In a nonzero product $p\,a_1\cdots a_k\,p$ with $a_1,\ldots,a_k\in pMp\cup\{h,h^*,u,u^*\}$, each $u^*=qu^*p$ is immediately preceded by $u=puq$, and $uu^*=p$. Hence $pM_{p,q}p=\Wstar(pMp,h)$. 
    
    Let $k\in\nn_{>0}$ and $n_1,\ldots,n_k\in\zz\setminus\{0\}$. Choose $x_0,\ldots,x_k\in pMp$ such that $\tau(x_i)=0$ for every $1\leq i\leq k-1$ and each of $x_0,x_k$ is either $p$ or has trace zero. Consider the alternating product
    $$x_0h^{n_1}x_1h^{n_2}\cdots h^{n_k}x_k.$$
    Substituting $h^n=(vu^*)^n$ and $h^{-n}=(uv^*)^n$ for $n\in\nn_{>0}$, we find $m\in\nn_{>0}$, $w_1,\ldots,w_m\in\{v,v^*\}$ and $y_0,\ldots,y_m\in M$ such that
    $$x_0h^{n_1}x_1h^{n_2}\cdots h^{n_k}x_k=y_0w_1y_1w_2\cdots w_my_m.$$
    Omit any initial or terminal factor equal to $p$. Each remaining $y_i$ is either $u$, $u^*$, or of the form $sx_jt$ with $s\in\{p,u^*\}$, $t\in\{p,u\}$ and $\tau(x_j)=0$. Since $\tau(u^*x_ju)=\tau(x_j)$, for every remaining $y_i$ we have
    $$y_i\in pMq\cup qMp\cup\{y\in pMp\cup qMq\mid\tau(y)=0\}.$$

    The remaining $y_i$ satisfy $\rE_D^M(y_i)=0$, while $\rE_D^{\rM_2(\cc)}(v)=\rE_D^{\rM_2(\cc)}(v^*)=0$. Therefore the remaining word is reduced in $M_{p,q}=M\ast_D\rM_2(\cc)$. By freeness with amalgamation over $D$, $\tau(y_0w_1y_1\cdots w_my_m)=0$. For a single power of $h$, this shows that $h$ is a Haar unitary, so $\Wstar(h)\cong\LL(\zz)$. In general, it shows that $h$ is free from $pMp$, which proves~\textup{(\ref*{it:equivalent-corners})}.
\end{proof}

{\emergencystretch=1em
Together, Lemma~\ref{lem:equivalent-endpoints} and Proposition~\ref{prop:scalar-amenability} leave only three possible amenable mergers, which the following proposition describes explicitly.\par}

\begin{Prop}
    \label{prop:amenable-merger}
    Let $(M,\tau)$ be a von Neumann algebra with a trace and let $p,q\in M$ be projections with $p+q=1$ and $0<\tau(p)=\tau(q)<\infty$. 
    
    Assume that $M_{p,q}$ is amenable. Then, up to exchanging the roles of $p$ and $q$, one of the following holds:
    \begin{enumerate}[\rm (i)]
        \item $pMq=0$, $pMp=\cc p$ and $qMq$ is amenable.  In this case, $qM_{p,q}q=qMq$ and $(M_{p,q},\tau)\cong\bigl(\rM_2(\cc)\ovt qMq,\ \Tr_2\pt\tau_{\mid qMq}\bigr).$
        \item $M\cong\cc^4$, and $p=h_1+h_2$, $q=h_3+h_4$ are sums of minimal projections such that $\tau(h_1)=\max_i\tau(h_i)$. In this case, $$\displaystyle (M_{p,q},\tau)\cong\underbrace{\rM_2(\cc)}_{\mathclap{(\tau(h_1)-\tau(h_4))\Tr_2}}\oplus\underbrace{\bigl(\rM_4(\cc)\ovt\rL^\infty([0,1],\nu)\bigr)}_{\tau(h_2)(\Tr_4\pt\nu)}\oplus\underbrace{\rM_2(\cc)}_{\mathclap{(\tau(h_1)-\tau(h_3))\Tr_2}}.$$ Moreover, each $h_i$ has rank one almost everywhere on the summand $\rM_4(\cc)\ovt\rL^\infty([0,1],\nu)$.
        \item $M\cong\rM_2(\cc)$, and $p,q$ are rank-one projections. In this case, $$(M_{p,q},\tau)\cong\bigl(\rM_2(\cc)\ovt\rL^\infty([0,1],\nu),\ \tau(p)(\Tr_2\pt\nu)\bigr).$$
    \end{enumerate}
\end{Prop}

\begin{proof}
    Since $M$ is an expected subalgebra of $M_{p,q}$, it is amenable. If $pMq=0$, then $M=pMp\oplus qMq$, and by Lemma~\ref{lem:equivalent-endpoints}\textup{(\ref*{it:separated-corners})} the algebra $pM_{p,q}p\cong(pMp)\ast(qMq)$ is amenable. By Proposition~\ref{prop:scalar-amenability}, either one of the two corners is scalar, and \textup{(i)} holds, or both corners are isomorphic to $\cc^2$, so that $M\cong\cc^4$ and \textup{(ii)} holds. In case~\textup{(ii)}, after exchanging $p$ and $q$ if necessary, label the minimal projections so that $\tau(h_1)=\max_i\tau(h_i)$.

    Suppose that $pMq\ne0$, and choose a nonzero partial isometry $u\in pMq$. We claim that $uu^*=p$. Otherwise $u^*u\ne q$ as well, because $\tau(uu^*)=\tau(u^*u)$ and $\tau(p)=\tau(q)$. The subalgebra $F=\Wstar(p,q,u)$ of $M$ is then isomorphic to $\rM_2(\cc)\oplus\cc\oplus\cc$: the matrix summand has matrix units $uu^*,u,u^*,u^*u$, and the scalar summands have units $p-uu^*$ and $q-u^*u$. In particular, neither $p$ nor $q$ is central in $F$. After normalizing the trace by $\tau(1_M)$, so that $\tau(p)=\tau(q)=1/2$, \cite[Lemma~4.2]{Dy95} (see also \cite[Lemma~2.5]{Re13}) shows that $F_{p,q}=F\ast_{\cc p\oplus\cc q}\rM_2(\cc)$ is the direct sum of a finite-dimensional algebra and of an algebra which is either an interpolated free group factor or a diffuse hyperfinite von Neumann algebra, the latter being possible only if $\dim_{\cc}F=4$. Since $\dim_{\cc}F=6$, the algebra $F_{p,q}$ has an interpolated free group factor as a direct summand, and is therefore not amenable. This contradicts its being an expected subalgebra of $M_{p,q}$ by Fact~\ref{fact:afp-subalgebras}, and proves the claim.

    Thus $uu^*=p$ and $u^*u=q$. Lemma~\ref{lem:equivalent-endpoints}\textup{(\ref*{it:equivalent-corners})} identifies $pM_{p,q}p$ with $(pMp)\ast\LL(\zz)$. Since this algebra is amenable and $\LL(\zz)$ is diffuse, Proposition~\ref{prop:scalar-amenability} forces $pMp=\cc p$. The matrix units $p,u,u^*,q$ now give $M\cong\rM_2(\cc)$, which is \textup{(iii)}.

For the displayed descriptions, Lemma~\ref{lem:equivalent-endpoints} gives $M_{p,q}\cong\rM_2(\cc)\ovt pM_{p,q}p$, with trace $\Tr_2\pt\tau_{\mid pM_{p,q}p}$. In case~\textup{(i)}, we have $pM_{p,q}p=v(qMq)v^*$ by the proof of Lemma~\ref{lem:equivalent-endpoints}\textup{(\ref*{it:separated-corners})}, so that $qM_{p,q}q=v^*(pM_{p,q}p)v=qMq$ and $pM_{p,q}p\cong qMq$. In case~\textup{(ii)}, apply Proposition~\ref{prop:two-projections} with $(\alpha,\beta,\gamma)=(\tau(h_1),\tau(h_3),\tau(p))$ and use $\tau(p)=\tau(h_1)+\tau(h_2){}=\tau(h_3)+\tau(h_4)$. Its last assertion, together with $h_3\sim vh_3v^*$ and $h_4\sim vh_4v^*$, shows that each $h_i$ has rank one almost everywhere on the diffuse summand. In case~\textup{(iii)}, $pM_{p,q}p\cong\LL(\zz)\cong\rL^\infty([0,1],\nu)$, with trace $\tau(p)\nu$. Amplifying gives the asserted trace-preserving isomorphisms, using $\Tr_2\pt\Tr_2=\Tr_4$ in case~\textup{(ii)}.
\end{proof}

As the main consequence of this proposition, we obtain that an amenable merger of a type~$\I$ algebra is again of type~$\I$.

\begin{Cor}
    \label{cor:amenable-merger-type-i}
    Let $(M,\tau)$ be a type~$\I$ von Neumann algebra with a trace, and let $p,q\in M$ be orthogonal projections with $0<\tau(p)=\tau(q)<\infty$. If $M_{p,q}$ is amenable, then $M_{p,q}$ is of type~$\I$.
\end{Cor}

\begin{proof}
    Put $r=p+q$ and $N=M_{p,q}$. Since $N$ is generated by $M$ and the adjoined matrix unit $v=rvr$, the projection $\csupp[M]{r}$ is central in $N$, so that
    $$
        N=\csupp[M]{r}N\oplus(1-\csupp[M]{r})M.
    $$
    By Lemma~\ref{lem:afp-identities}\textup{(\ref*{it:afp-compression})}, $rNr\cong(rMr)_{p,q}$, a merger in the finite tracial algebra $rMr$ with $p+q=1_{rMr}$. The corner $rNr$ is amenable by Fact~\ref{fact:amenability-properties}\textup{(\ref*{it:amenability-expected-corners})}, so Proposition~\ref{prop:amenable-merger} identifies it with one of the three displayed algebras. Each of them is of type~$\I$, since in case~\textup{(i)} the algebra $qMq$ is a corner of the type~$\I$ algebra $M$. Since $\csupp[N]{r}$ dominates every $M$-unitary conjugate of $r$, we have $\csupp[N]{r}=\csupp[M]{r}$, so $r$ is full in $\csupp[M]{r}N$, which is therefore of type~$\I$. Hence $N$ is of type~$\I$.
\end{proof}

\subsection{Controlling traces under merging}
\label{subsec:afp-mergers}

In the proof of Theorem~\ref{thm:lattice}, the amalgamated free product will be obtained as the weak$^\star$ closure of an increasing union of successive mergers (Subsection~\ref{subsec:lattice-proof}). By Corollary~\ref{cor:amenable-merger-type-i}, each amenable merger of a type~$\I$ algebra is of type~$\I$. To prove that this closure is of type~$\I$, we must also exclude a hyperfinite type~$\II$ summand in the limit, since the weak$^\star$ closure of an increasing union of type~$\I$ algebras can be of type~$\II$. We do so by uniformly controlling the traces of enough projections at each step. This is the purpose of the next definitions, which single out projections whose traces are not too small, in a way that is preserved by mergers.

\begin{Def}
    \label{def:integral-atoms}
    Let $(M,\tau)$ be a type~$\I$ von Neumann algebra with a trace. We say that: 
    \begin{enumerate}[\rm (i)]
        \item\label{it:integral-atoms} $M$ \emph{has integral atoms} if $\tau(f)\in\nn$ for every minimal projection $f\in M$.
        \item\label{it:rank-controlled} A finite-trace projection $e\in M$ is \emph{rank-controlled in $M$} if $\rk_M(e\zdiff[M])\leq\tau(e\zdiff[M])$ (equivalently, if $(e\zdiff[M])M(e\zdiff[M])$ is of type $\I_{\leq\tau(e\zdiff[M])}$).
    \end{enumerate}
\end{Def}

The condition on minimal projections is stricter than a more symmetric assumption such as ``$\rk(f)\leq\tau(f)$ for every atomic projection $f\in M$''. The reason is that free product formulas produce differences of traces, as witnessed by Proposition~\ref{prop:two-projections}, so such a property would not be preserved by mergers. The next lemma records how these notions behave under corners, Murray--von Neumann equivalence and orthogonal sums. 

\begin{Lem}
    \label{lem:rank-control-properties}
    Let $(M,\tau)$ be a type $\I$ von Neumann algebra with a trace. 
    \begin{enumerate}[\rm (i)]
        \item\label{it:integral-atoms-corners} Let $p\in M$ be a projection. Then $pMp$ has integral atoms if and only if $\csupp{p}M$ has integral atoms.
        \item\label{it:rank-control-corners} Let $p\in M$ be a finite-trace projection. Then $p$ is rank-controlled in $M$ if and only if it is rank-controlled in $pMp$.
        \item\label{it:rank-control-equivalence} Let $p,q\in M$ be two Murray--von Neumann equivalent finite-trace projections. Then $p$ is rank-controlled in $M$ if and only if $q$ is rank-controlled in $M$. 
        \item\label{it:rank-control-sums}  Let $p,q\in M$ be two orthogonal rank-controlled finite-trace projections. Then $p+q$ is rank-controlled.
    \end{enumerate}
\end{Lem}

\begin{proof}
    \textup{(\ref*{it:integral-atoms-corners})} Every minimal projection $e\in pMp$ is minimal in $\csupp{p}M$, since its corner in either algebra is $eMe$. Conversely, let $e\in\csupp{p}M$ be a minimal projection. Since $\csupp[M]{e}\leq\csupp{p}$, we have $eMp\ne0$. The polar decomposition of a nonzero element of $eMp$ gives a partial isometry $v\in M$ with $vv^*=e$ and $v^*v\leq p$. The projection $v^*v$ is minimal in $pMp$ and satisfies $\tau(v^*v)=\tau(e)$. Thus the traces of minimal projections of $pMp$ and of $\csupp{p}M$ form the same set, which proves the equivalence.

    \textup{(\ref*{it:rank-control-corners})} The atomic and diffuse parts of $pMp$ are the corresponding corners of those of $M$, so $\zdiff[pMp]=p\zdiff[M]$. The diffuse corner entering either rank-control condition is therefore $(p\zdiff[M])M(p\zdiff[M])$, with the same restricted trace.

    \textup{(\ref*{it:rank-control-equivalence})} A partial isometry implementing the equivalence between $p$ and $q$ induces a trace-preserving isomorphism $pMp\cong qMq$, so the assertion follows from~\textup{(\ref*{it:rank-control-corners})}.

    \textup{(\ref*{it:rank-control-sums})} We have $\tau(p+q)<\infty$, and additivity of rank and trace gives
    \begin{align*}
            \rk_M((p+q)\zdiff[M])
            &=\rk_M(p\zdiff[M])+\rk_M(q\zdiff[M])\\
            &\leq\tau(p\zdiff[M])+\tau(q\zdiff[M])
            =\tau((p+q)\zdiff[M]).\qedhere
    \end{align*}
\end{proof}

These notions are designed to yield the following uniform bound on matrix sizes, which passes to weak$^\star$ limits by Lemma~\ref{lem:standard-identity}.

\begin{Lem}
    \label{lem:rank-controlled}
    Let $(M,\tau)$ be a type $\I$ von Neumann algebra with a trace. Suppose that $M$ has integral atoms. If $e\in M$ is a rank-controlled projection, then $eMe$ is of type $\I_{\leq\tau(e)}$.
\end{Lem}

\begin{proof}
    By Lemma~\ref{lem:rank-control-properties}\textup{(\ref*{it:integral-atoms-corners})--(\ref*{it:rank-control-corners})}, $eMe$ has integral atoms and its unit $e$ is rank-controlled. Its diffuse part is therefore of type $\I_{\leq\tau(e)}$. Since $\tau(e)<\infty$, every atomic factor of $eMe$ is a matrix algebra, whose size is at most $\tau(e)$ because each minimal projection has trace at least $1$. Thus $eMe$ is of type $\I_{\leq\tau(e)}$.
\end{proof}

We now prove a quantitative version of Corollary~\ref{cor:amenable-merger-type-i}.

\begin{Prop}
    \label{prop:amenable-free-merger}
    Let $(M,\tau)$ be a type~$\I$ von Neumann algebra with a trace.  Let $p,q\in M$ be orthogonal projections such that $0<\tau(p)=\tau(q)<\infty$. Assume that $(M,\tau)$ has integral atoms.

    If $M_{p,q}$ is amenable, then $(M_{p,q},\tau)$ is of type~$\I$ and has integral atoms. Moreover, if $p$ and $q$ are rank-controlled in $M$, then every rank-controlled projection of $M$ remains rank-controlled in $M_{p,q}$.
\end{Prop}

\begin{proof}
 Set $N=M_{p,q}$, assume that $N$ is amenable and set $r=p+q$. By Fact~\ref{fact:amenability-properties}\textup{(\ref*{it:amenability-expected-corners})}, $rNr$ is amenable. As in the proof of Corollary~\ref{cor:amenable-merger-type-i}, $\csupp[M]{r}=\csupp[N]{r}=\csupp{r}$ and
    \begin{equation}\label{eq:merger-central-decomposition}
        N=\csupp{r}N\oplus(1-\csupp{r})M.
    \end{equation}
    Only the central summand $\csupp{r}N$ therefore differs from $M$.

    \smallskip
    \noindent\textbf{Integral atoms.}\enspace By Corollary~\ref{cor:amenable-merger-type-i}, $N$ is of type~$\I$. Lemma~\ref{lem:afp-identities}\textup{(\ref*{it:afp-compression})} gives $rNr\cong rMr\ast_{\cc p\oplus\cc q}\rM_2(\cc)$, so Proposition~\ref{prop:amenable-merger} applies to $rMr$. Its descriptions show that $rNr$ has integral atoms: the atomic traces are inherited from a corner of $M$ in case~\textup{(i)}, are integer differences of the $\tau(h_i)$ in case~\textup{(ii)}, and there are no atoms in case~\textup{(iii)}. Since $r$ is full in $\csupp{r}N$, Lemma~\ref{lem:rank-control-properties}\textup{(\ref*{it:integral-atoms-corners})} shows that $\csupp{r}N$ has integral atoms. As $(1-\csupp{r})N=(1-\csupp{r})M$, the algebra $N$ has integral atoms.

    \smallskip
    \noindent\textbf{Rank control.}\enspace We establish the following claims:
    \begin{claims}
        \label{fact:merger-key-facts}
        \leavevmode
        \begin{enumerate}[\rm (a)]
            \item\label{it:merger-fact-diffuse} $\zdiff[M]N\zdiff[M]=\zdiff[M]M$.
            \item\label{it:merger-fact-rank} If $p$ and $q$ are rank-controlled in $M$, then every minimal projection of $M$ is rank-controlled in $N$.
        \end{enumerate}
    \end{claims}
    This gives the conclusion: assume that $p$ and $q$ are rank-controlled in $M$ and let $e\in M$ be a rank-controlled projection. Decompose $e=e\zdiff[M]+e\zat[M]$. By Lemma~\ref{lem:rank-control-properties}\textup{(\ref*{it:rank-control-sums})}, it suffices to prove that $e\zdiff[M]$ and $e\zat[M]$ are rank-controlled in $N$. By Lemma~\ref{lem:rank-control-properties}\textup{(\ref*{it:rank-control-corners})} and~\textup{(\ref*{it:merger-fact-diffuse})}, $e\zdiff[M]$ is rank-controlled in $N$. As $\tau(e)<\infty$ and $M$ has integral atoms, $e\zat[M]$ is a finite orthogonal sum of minimal projections of $M$. By~\textup{(\ref*{it:merger-fact-rank})} and Lemma~\ref{lem:rank-control-properties}\textup{(\ref*{it:rank-control-sums})}, it is rank-controlled in $N$.

    We now prove claims~\textup{(\ref*{it:merger-fact-diffuse})} and~\textup{(\ref*{it:merger-fact-rank})} in the three cases of Proposition~\ref{prop:amenable-merger} applied to $rMr$.

    \begin{enumerate}[\rm (i)]
    \item Assume, as in Proposition~\ref{prop:amenable-merger}\textup{(i)}, that $pMq=0$ and $pMp=\cc p$. By Proposition~\ref{prop:amenable-merger}\textup{(i)}, applied to $rMr$, we have $qNq=qMq$. As $q$ is full in $\csupp[M]{q}M$, the equality $qNq=qMq$ implies $\csupp[M]{q}N\csupp[M]{q}=\csupp[M]{q}M$.  Since $pMq=0$, we have $\csupp[M]{p}\csupp[M]{q}=0$, hence $\csupp{r}=\csupp[M]{p}+\csupp[M]{q}$. Together with \eqref{eq:merger-central-decomposition}, this gives $$(1-\csupp[M]{p})N(1-\csupp[M]{p})=(1-\csupp[M]{p})M.$$ As $p$ is minimal in $M$, $\csupp[M]{p}\leq\zat[M]$, so that $\zdiff[M]\leq1-\csupp[M]{p}$. Therefore, $\zdiff[M]N\zdiff[M]=\zdiff[M]M$. 
    
    Every minimal projection of $(1-\csupp[M]{p})M$ remains minimal in $N$, hence is rank-controlled. Every minimal projection of $\csupp[M]{p}M$ is equivalent to $p$ in $M$, hence to $q$ in $N$. Since $qNq=qMq$ and $q$ is rank-controlled in $M$, Lemma~\ref{lem:rank-control-properties}\textup{(\ref*{it:rank-control-corners})--(\ref*{it:rank-control-equivalence})} shows that these projections are rank-controlled in $N$.

    \item Assume that $rMr\cong\cc^4$, and use the notation of Proposition~\ref{prop:amenable-merger}\textup{(ii)}. As $M$ has integral atoms, $\tau(h_i)\in \nn_{>0}$ for $i\in \{1,\dots, 4\}$. Since $h_iMh_j\subset rMr$, the central supports of the $h_i$ in $M$ are pairwise orthogonal. Therefore $\csupp{r}=\sum_j\csupp[M]{h_j}\leq\zat[M]$ is an atomic projection. Equation \eqref{eq:merger-central-decomposition} gives $\zdiff[M]N\zdiff[M]=\zdiff[M]M$ and shows that every minimal projection of $(1-\csupp{r})M$ remains minimal in $N$, hence rank-controlled. 
    
    Every minimal projection of $\csupp{r}M$ is equivalent in $M$ to some $h_j$. By Proposition~\ref{prop:amenable-merger}\textup{(ii)}, applied to $rMr$, each $h_i$ has rank one almost everywhere on the diffuse summand of $rNr$, where its trace is $\tau(h_2)\geq1$. Hence the $h_i$ are rank-controlled in $rNr$, and therefore in $N$ by Lemma~\ref{lem:rank-control-properties}\textup{(\ref*{it:rank-control-corners})}. By Lemma~\ref{lem:rank-control-properties}\textup{(\ref*{it:rank-control-equivalence})}, the same holds for every minimal projection of $\csupp{r}M$.

    \item Assume that $rMr\cong\rM_2(\cc)$. By Proposition~\ref{prop:amenable-merger}\textup{(iii)}, applied to $rMr$, and since $p\sim q$ in $N$, the corner $pNp$ is abelian and diffuse.
    As $p$ and $q$ are equivalent minimal projections of $M$, $\csupp{r}M$ is a type~$\I$ factor and $\csupp{r}\leq\zat[M]$. Equation \eqref{eq:merger-central-decomposition} then gives
    $\zdiff[M]N\zdiff[M]=\zdiff[M]M$.

    Every minimal projection of $(1-\csupp{r})M$ remains minimal in $N$, hence rank-controlled. Every minimal projection of $\csupp{r}M$ is equivalent to $p$. Since $pNp$ is abelian and diffuse, $\rk_{pNp}(p)=1\leq\tau(p)$ on its diffuse part, so $p$ is rank-controlled in $pNp$. Lemma~\ref{lem:rank-control-properties}\textup{(\ref*{it:rank-control-corners})--(\ref*{it:rank-control-equivalence})} therefore shows that every minimal projection of $\csupp{r}M$ is rank-controlled in $N$.\qedhere
    \end{enumerate}
\end{proof}

\subsection{Proof of Theorem~\ref{thm:lattice}}
\label{subsec:lattice-proof}

We first treat the case of abelian amalgams in the next proposition. 

\begin{Prop}
    \label{prop:abelian-amalgam}
    Let $(A_1\supset D\subset A_2,\tau)$ be tracial inclusions of atomic von Neumann algebras. Assume that $D$ is abelian and that $\tau(e)\in\nn$ for every minimal projection $e\in A_i$, $i\in\{1,2\}$. 
    
    If $M=A_1\ast_D A_2$ is amenable, then $M$ is of type~$\I$, and $eMe$ is of type $\I_{\leq\tau(e)}$ for every minimal projection $e\in A_i$, $i\in\{1,2\}$.
\end{Prop}

\begin{proof}
    For each $i\in\{1,2\}$, fix an atomic masa $D_i\subset A_i$ containing $D$. 
    In each block of $A_i$, choose matrix units $e_{jl}\in A_i$ whose diagonal projections are the minimal projections of $D_i$ in this block. Enumerate the units $e_{j1}$ with $j\in\nn_{\geq2}$, over all blocks of $A_i$, as a finite or infinite sequence $(e^i_k)_k$. Since $e_{jl}=e_{j1}e_{l1}^*$, we have $A_i=D_i\vee\bigvee_k e^i_k$.

    Enumerate the two families $(e^1_k)_k$ and $(e^2_k)_k$ together, without repetitions, as a single sequence $(v_n)$ in $A_1\cup A_2$, indexed from $1$ and either finite or infinite. For $n\in\nn$, let $A_{i,n}$ be the von Neumann algebra generated by $D_i$ and those $v_k$ with $k\leq n$ chosen from $A_i$, and put
    $$
        N_n=A_{1,n}\ast_D A_{2,n},\qquad\text{so that}\qquad
        N_0=D_1\ast_D D_2.
    $$
    By Fact~\ref{fact:afp-subalgebras}, the inclusions $N_n\subset N_{n+1}\subset M$ are trace-preserving, hence expected, so every $N_n$ is amenable by Fact~\ref{fact:amenability-properties}\textup{(\ref*{it:amenability-expected-corners})}. Since $A_i$ is generated by $\bigcup_nA_{i,n}$, we have
    \begin{equation}\label{eq:dense-exhaustion}
        M=\Wstar\Bigl(\bigcup_{n\geq0}N_n\Bigr).
    \end{equation}

    \smallskip
    \noindent\textbf{Step 1.}\enspace We show that $N_0$ is of type~$\I$, has integral atoms, and that every minimal projection of $D_1$ and $D_2$ is rank-controlled in $N_0$. Since $D_1$ and $D_2$ are abelian, $D$ is central in $N_0$. Let $z\in D$ be a minimal projection. Then $\tau(z)<\infty$, as $\tau$ is semifinite on $D$, and compressing reduced words by $z$ shows that
    $$
        zN_0\cong zD_1\ast zD_2.
    $$
    The minimal projections of $zD_1$ and $zD_2$ are minimal projections of $A_1$ and $A_2$, so their traces are positive integers. In particular, $zD_1$ and $zD_2$ are finite dimensional. As $zN_0$ is amenable, Proposition~\ref{prop:scalar-amenability} shows that one of $zD_1$ and $zD_2$ equals $\cc z$, or that both are isomorphic to $\cc^2$.

    If $zD_2=\cc z$, say, then $zN_0=zD_1$ is abelian and finite dimensional, and its minimal projections have traces in $\nn$. It has no diffuse part, so every finite-trace projection of $zN_0$ is rank-controlled in $zN_0$.

    Otherwise, let $h_1,h_2\in zD_1$ and $h_3,h_4\in zD_2$ be their respective minimal projections. Interchanging $D_1$ and $D_2$ and relabeling the $h_j$ if necessary, we may assume that $\tau(h_1)=\max_j\tau(h_j)$. Proposition~\ref{prop:two-projections} then gives
    $$
        \begin{aligned}
            zN_0&\cong\left(\underset{\tau(h_1)}{\cc}\oplus\underset{\tau(h_2)}{\cc}\right)
            \ast\left(\underset{\tau(h_3)}{\cc}\oplus\underset{\tau(h_4)}{\cc}\right)\\
            &\cong\underset{\tau(h_1)-\tau(h_4)}{\cc}
            \oplus\underbrace{\bigl(\rM_2(\cc)\ovt\rL^\infty([0,1],\nu)\bigr)}_{\tau(h_2)(\Tr_2\pt\nu)}
            \oplus\underset{\tau(h_1)-\tau(h_3)}{\cc}.
        \end{aligned}
    $$
    The traces of the minimal projections of $zN_0$ are therefore integers. By the last assertion of Proposition~\ref{prop:two-projections}, each $h_j$ has rank one almost everywhere on the diffuse part of $zN_0$, where its trace is $\tau(h_2)\geq1$. Hence each $h_j$ is rank-controlled in $zN_0$.

    In both cases, $zN_0$ is of type~$\I$ and has integral atoms, and the minimal projections of $D_1$ and $D_2$ below $z$ are rank-controlled in $zN_0$. Since the minimal projections $z\in D$ are central in $N_0$ with sum $1$, the algebra $N_0=\bigoplus_z zN_0$ is of type~$\I$ and has integral atoms. Every minimal projection of $D_1$ or $D_2$ lies below some $z$, so it is rank-controlled in $N_0$ by Lemma~\ref{lem:rank-control-properties}\textup{(\ref*{it:rank-control-corners})}.

    \smallskip
    \noindent\textbf{Step 2.}\enspace We show by induction on $n$ that $N_n$ is of type~$\I$, has integral atoms, and that every minimal projection of $D_1$ and $D_2$ is rank-controlled in $N_n$. The case $n=0$ is Step~1.

    Assume that this holds for $N_n$. If $v_{n+1}$ is not defined, then $N_{n+1}=N_n$. Otherwise, $v_{n+1}=e_{j1}$ for some $j\in\nn_{\geq2}$, where $(e_{jl})$ is the system of matrix units chosen in a block of $A_i$. The projections $e_{jj}$ and $e_{11}$ are minimal projections of $D_i$ with the same trace. Since $e_{j1}$ has not been adjoined before, $e_{jj}$ supports a scalar block of $A_{i,n}$, distinct from the block containing $e_{11}$. 
    Example~\ref{ex:block-merger} therefore identifies $A_{i,n+1}$ with $(A_{i,n})_{e_{jj},e_{11}}$ and Lemma~\ref{lem:outer-fusion} identifies $N_{n+1}$ with the merger $(N_n)_{e_{jj},e_{11}}$. By the induction hypothesis, $e_{jj}$ and $e_{11}$ are rank-controlled in $N_n$. Since $N_{n+1}$ is amenable, Proposition~\ref{prop:amenable-free-merger} therefore shows that $N_{n+1}$ is of type~$\I$, has integral atoms, and that every minimal projection of $D_1$ and $D_2$ remains rank-controlled in $N_{n+1}$.

    \smallskip
    \noindent\textbf{Step 3.}\enspace Let $i\in\{1,2\}$ and let $e\in D_i$ be a minimal projection, so that $\tau(e)\in\nn$. By Step~2 and Lemma~\ref{lem:rank-controlled}, $eN_ne$ is of type $\I_{\leq\tau(e)}$ for every $n\in\nn$. Lemma~\ref{lem:standard-identity} shows that $S_{2\tau(e)}$ vanishes identically on $\bigcup_neN_ne$. By \eqref{eq:dense-exhaustion}, this union is a weak$^\star$ dense $*$-subalgebra of $eMe$, and Lemma~\ref{lem:standard-identity} shows that $eMe$ is of type $\I_{\leq\tau(e)}$.

    The minimal projections $e\in D_1$ sum to $1$, so the central supports $\csupp[M]{e}$ have supremum $1$. Since $e$ is full in $\csupp[M]{e}M$ and $eMe$ is of type~$\I$, each $\csupp[M]{e}M$ is of type~$\I$, and so is $M$. Finally, every minimal projection of $A_i$ is equivalent in $A_i$ to a minimal projection of $D_i$, and equivalent projections have the same trace and isomorphic corners.
\end{proof}

\begin{proof}[Proof of Theorem~\ref{thm:lattice}]
    Assume that $M$ is amenable. Replacing $\tau$ by $\delta^{-1}\tau$ on $A_1$, $A_2$ and $B$ does not change the trace-preserving conditional expectations, hence does not change $M$, and multiplies the trace of $M$ by $\delta^{-1}$. We may therefore assume that $\delta=1$. We must show that $M$ is of type~$\I$, and that $pMp$ is of type $\I_{\leq\tau(p)}$ for every finite-trace projection $p\in A_1\cup A_2$.

    Let $f\in B$ be the sum of a choice of one minimal projection in each block of $B$. Then $fBf$ is abelian and $\csupp[B]{f}=1$, so Lemma~\ref{lem:afp-identities}\textup{(\ref*{it:afp-full-corner})} gives a trace-preserving isomorphism
    $$
        fMf\cong fA_1f\ast_{fBf}fA_2f.
    $$
    The algebras $fA_1f$ and $fA_2f$ are atomic, and their minimal projections are minimal projections of $A_1$ and $A_2$, with traces in $\nn$. By Fact~\ref{fact:amenability-properties}\textup{(\ref*{it:amenability-expected-corners})}, the corner $fMf$ is amenable. Proposition~\ref{prop:abelian-amalgam} therefore shows that $fMf$ is of type~$\I$, and that $eMe$ is of type $\I_{\leq\tau(e)}$ for every minimal projection $e\in fA_if$, $i\in\{1,2\}$. Since $\csupp[M]{f}\geq\csupp[B]{f}=1$, the projection $f$ is full in $M$, so $M$ is of type~$\I$.

    Let $p\in A_i$ be a finite-trace projection, for some $i\in\{1,2\}$. As the minimal projections of $A_i$ have trace at least $1$, $p$ is a finite orthogonal sum of minimal projections of $A_i$. Let $e\in A_i$ be one of them. Since $\csupp[A_i]{f}\geq\csupp[B]{f}=1$, the block of $A_i$ containing $e$ contains a minimal projection of $fA_if$, which is equivalent to $e$ in $A_i$. Hence $eMe$ is of type $\I_{\leq\tau(e)}$, that is, $\rk_M(e)\leq\tau(e)$ almost everywhere. By additivity of rank and trace, $\rk_M(p)\leq\tau(p)$ almost everywhere, so $pMp$ is of type $\I_{\leq\tau(p)}$. For the original trace, this is the bound $\tau(p)/\delta$.
\end{proof}

\section{From \texorpdfstring{$\Wstar$}{W*}-amenability to uniform admissibility}
\label{sec:compact-amalgams}
In this section, we prove that a $\Wstar$-amenable amalgam of compact groups over a common compact open subgroup is uniformly admissible, hence of type~$\I$.

\begin{The}
    \label{thm:compact-amalgam-ccr}
    Let $K_1$ and $K_2$ be second countable totally disconnected compact groups with a common open subgroup $J$. If $K_1\ast_JK_2$ is $\Wstar$-amenable, then it is uniformly admissible. More precisely, with $d_i=[K_i:J]\in\nn_{>0}$ for each $i\in\{1,2\}$, every irreducible unitary representation $\pi:K_1\ast_JK_2\to\U(\H_\pi)$ satisfies, for every compact open subgroup $U\leq J$,
    $$
        \dim\H_\pi^U\leq\lcm\left(d_1,d_2\right)[J:U].
    $$
\end{The}

\subsection{Compact amalgams and the Plancherel trace}
\label{subsec:compact-plancherel}

To apply Theorem~\ref{thm:lattice}, we need the traces of minimal projections in the three compact group algebras to lie in a common discrete subgroup of $\rr$. The following computation expresses these traces in terms of representation dimensions and subgroup indices.

\begin{Lem}
    \label{lem:compact-amalgam-afp}
    Let $K_1$ and $K_2$ be second countable compact groups with a common open subgroup $J$. Normalize Haar measure by $\mu(J)=1$. For $K\in\{K_1,K_2,J\}$ and $\pi\in\widehat K$, every minimal projection in the corresponding matrix block of $\LL(K)$ has Plancherel trace $\dim\pi/[K:J]$.
\end{Lem}

\begin{proof}
    Put $G=K_1\ast_JK_2$ and normalize Haar measure by $\mu_G(J)=1$.

    Let $K\in\{K_1,K_2,J\}$, with Haar measure $\mu_G|_K$. Since $K$ is a second countable compact group, its unitary dual is countable and every irreducible representation of $K$ is finite dimensional. Choose one representative of each class in $\widehat K$. By Peter--Weyl theory \cite[Theorem~5.12]{Fo16}, for every $k\in K$ we have
    $$
        \begin{gathered}
            \rL^2(K)\cong\bigoplus_{\sigma\in\widehat K}
            \H_\sigma\pt\ovl{\H_\sigma},\\
            \lambda_K(k)\cong\bigoplus_{\sigma\in\widehat K}
            \bigl(\sigma(k)\pt1_{\ovl{\H_\sigma}}\bigr).
        \end{gathered}
    $$
    The second tensor factor accounts for the multiplicity $\dim(\sigma)$ of $\sigma$. Consequently,
    $$
        \LL(K)\cong\bigoplus_{\sigma\in\widehat K}
        \bigl(\rB(\H_\sigma)\ovt\cc1_{\ovl{\H_\sigma}}\bigr)
        \cong\bigoplus_{\sigma\in\widehat K}\rB(\H_\sigma).
    $$
    The minimal projections of $\LL(K)$ correspond to rank-one projections in a single matrix block $\rB(\H_\sigma)$.

    Fix $\pi\in\widehat K$ and a unit vector $\eta\in\H_\pi$, and let $q$ be the orthogonal projection onto $\cc\eta$. We use inner products linear in the first variable, so $q\xi=\ps{\xi}{\eta}\eta$. For every $k\in K$, define
    $$
        f(k)=\frac{\dim(\pi)}{\mu_G(K)}\ovl{\ps{\pi(k)\eta}{\eta}},
    $$
    and extend $f$ by zero outside $K$. Since $K$ is compact and open, $f\in\rC_c(G)$ and $\lambda_G(f)\in\LL(K)\subset\LL(G)$. Integrating the formula for $\lambda_K(k)$ above shows that the component of $\lambda_G(f)$ in the block corresponding to $\sigma$ is
    $$
        \left(\int_K f(k)\sigma(k)\,{\rm d}\mu_G(k)\right)
        \pt1_{\ovl{\H_\sigma}}.
    $$
    For $\xi,\zeta\in\H_\sigma$, Schur orthogonality \cite[Theorem~5.8]{Fo16} gives
    $$
        \int_K\ovl{\ps{\pi(k)\eta}{\eta}}\ps{\sigma(k)\xi}{\zeta}\,{\rm d}\mu_G(k)
        =\begin{cases}
            \dfrac{\mu_G(K)}{\dim(\pi)}\ps{\xi}{\eta}\ps{\eta}{\zeta} & \text{if }\sigma=\pi,\\[6pt]
            0 & \text{if }\sigma\ne\pi.
        \end{cases}
    $$
    Multiplying by $\dim(\pi)/\mu_G(K)$ and comparing matrix coefficients yields
    $$
        \int_K f(k)\sigma(k)\,{\rm d}\mu_G(k)
        =\begin{cases}
            q & \text{if }\sigma=\pi,\\
            0 & \text{if }\sigma\ne\pi.
        \end{cases}
    $$
    Hence $\lambda_G(f)$ equals $q\pt1_{\ovl{\H_\pi}}$ in the block corresponding to $\pi$ and vanishes in every other block. As $\pi$ and $\eta$ vary, these are all the minimal projections of $\LL(K)$.

    Taking $\sigma=\pi$ and $\xi=\zeta=\eta$ in the Schur orthogonality formula gives
    $$
        \int_K\abs{\ps{\pi(k)\eta}{\eta}}^2\,{\rm d}\mu_G(k)
        =\frac{\mu_G(K)}{\dim(\pi)}.
    $$
    Since $\lambda_G(f)$ is a projection, the Plancherel trace identity yields
    $$
        \begin{aligned}
            \tau_G(\lambda_G(f))
            &=\tau_G\bigl(\lambda_G(f)^*\lambda_G(f)\bigr)=\norm{f}_2^2\\
            &=\frac{\dim(\pi)^2}{\mu_G(K)^2}\cdot\frac{\mu_G(K)}{\dim(\pi)}
            =\frac{\dim(\pi)}{\mu_G(K)}.
        \end{aligned}
    $$
    Since $\mu_G(K)=[K:J]$, the asserted formula follows. In particular, all these traces belong to $\delta\nn$, where $\delta=1/\lcm([K_1:J],[K_2:J])$.
\end{proof}

\subsection{From corner bounds to uniform admissibility}
\label{subsec:compact-admissibility}

Put $G=K_1\ast_JK_2$. If $\LL(G)$ is amenable, then Theorem~\ref{thm:lattice} gives a bound on the matrix sizes in $p_U\LL(G)p_U$ for every compact open subgroup $U\leq J$. Hence the Hecke algebra $\Hk GU$ satisfies a standard polynomial identity, which passes to its image under every irreducible unitary representation of $G$.

\begin{Prop}
    \label{prop:full-pi-transfer}
    Let $G$ be a second countable totally disconnected locally compact group. Suppose that there is a neighborhood basis $\U$ of the identity consisting of compact open subgroups such that, for every $U\in\U$, the corner $p_U\LL(G)p_U$ is of type $\I_{\leq n_U}$ for some $n_U\in\nn$. Then every irreducible unitary representation $\pi:G\to\U(\H_\pi)$ satisfies, for every $U\in\U$,
    $$
        \dim\H_\pi^U\leq n_U.
    $$
    Consequently, $G$ is uniformly admissible.
\end{Prop}

\begin{proof}
    Fix $U\in\U$. Using the injectivity of $\lambda_G$ recalled in Subsection~\ref{subsec:prelim-group-vna}, we view $\Hk GU$ as a subalgebra of $p_U\LL(G)p_U$. For all $f_1,\ldots,f_{2n_U}\in\Hk GU$, Lemma~\ref{lem:standard-identity} then gives
    $$
        S_{2n_U}(f_1,\ldots,f_{2n_U})=0.
    $$
    Here the polynomial is evaluated using convolution as multiplication. Let $\pi:G\to\U(\H_\pi)$ be an irreducible unitary representation and put $E=\pi(f_U)\in\rB(\H_\pi)$, whose range is $\H_\pi^U$. Applying $\pi$ to the convolution identity shows that $S_{2n_U}$ vanishes on $\pi(\Hk GU)$. If $E=0$, the required dimension bound is immediate. Suppose that $E\ne0$. As in Subsection~\ref{subsec:prelim-group-vna}, with $\pi$ in place of $\lambda_G$, irreducibility implies that $\pi(\Hk GU)$ is weak$^\star$ dense in $\rB(E\H_\pi)$. Lemma~\ref{lem:standard-identity} therefore gives $\dim\H_\pi^U\leq n_U$. No weak containment assumption on $\pi$ has been used.

    For an arbitrary compact open subgroup $V\leq G$, choose $U\in\U$ with $U\leq V$. Then $\H_\pi^V\subset\H_\pi^U$, which proves uniform admissibility.
\end{proof}

We can now combine the Plancherel trace computation with the discrete trace criterion and the preceding proposition to prove uniform admissibility for compact amalgams.

\begin{proof}[Proof of Theorem~\ref{thm:compact-amalgam-ccr}]
    Put $G=K_1\ast_JK_2$, which is a second countable totally disconnected locally compact group by Subsection~\ref{subsec:prelim-bass-serre}. Since $G$ is $\Wstar$-amenable, $\LL(G)$ is amenable. Normalize the Haar measure by $\mu_G(J)=1$. By \eqref{eq:compact-amalgam-afp} and Lemma~\ref{lem:compact-amalgam-afp}, $\LL(G)$ is the amalgamated free product of the atomic algebras $\LL(K_1)$ and $\LL(K_2)$ over $\LL(J)$, and \eqref{eq:lattice-hypothesis} holds with $\delta=1/\lcm(d_1,d_2)$. These inclusions are expected, since the Plancherel trace is semifinite on $\LL(J)$. If $U\leq J$ is compact open, then $p_U\in\LL(J)$ and $\tau_G(p_U)=\mu_G(U)^{-1}=[J:U]$. Theorem~\ref{thm:lattice} shows that $p_U\LL(G)p_U$ is of type $\I_{\leq\lcm(d_1,d_2)[J:U]}$. Since the compact open subgroups of $J$ form a neighborhood basis of the identity, Proposition~\ref{prop:full-pi-transfer} proves the result.
\end{proof}

In the applications below, the compact amalgam description applies to the type-preserving subgroup. The following fact allows us to pass from this subgroup to the whole group.

\begin{fact}
    \label{fact:uniform-finite-extension}
    Let $G$ be a second countable totally disconnected locally compact group and $H\leq G$ an open subgroup of index at most two. If $H$ is uniformly admissible, then so is $G$.
\end{fact}

This is \cite[Lemma~A.8]{Semal24}.

\subsection{Proof of Theorem~\ref{thm:main}}
\label{subsec:main-proofs}

\begin{proof}
    We apply Lemma~\ref{lem:tree-thick-reduction}. In cases~\textup{(\ref*{it:tree-finite})} and~\textup{(\ref*{it:tree-line})}, $G$ is discrete and virtually abelian. Its irreducible representations have uniformly bounded finite dimensions, so $G$ is uniformly admissible.

    In case~\textup{(\ref*{it:tree-thick})}, we identify $G$ with a closed subgroup of $\Aut(S)$, whose action on $\partial S$ is two-transitive. Let $G^+\leq G$ be the type-preserving subgroup. By \cite[Lemma~3.1.1]{BM00}, $G$ is locally transitive and its vertex stabilizers act transitively on spheres. Choose adjacent $v_1,v_2\in V(S)$ and put $K_i=(G^+)_{v_i}$ for each $i\in\{1,2\}$ and $J=K_1\cap K_2$. By Lemma~\ref{lem:tree-amalgam}, $G^+=K_1\ast_JK_2$ is an amalgam of compact open subgroups. Its vertex stabilizers coincide with those of $G$ and hence act transitively on spheres. In particular, a sphere about $v_1$ that meets the orbit $G^+v_1$ is contained in it, and $K_1$ acts transitively on this sphere. Proposition~\ref{prop:ttn-spherical-amenability}, applied to $G^+$, shows that $G^+$ is $\Wstar$-amenable. Theorem~\ref{thm:compact-amalgam-ccr} shows that $G^+$ is uniformly admissible, and Fact~\ref{fact:uniform-finite-extension} shows that $G$ is uniformly admissible.
\end{proof}

\subsection{Proof of Theorem~\ref{thm:equivalences}}
\label{subsec:equivalences-proof}

\begin{proof}
    Theorem~\ref{thm:main} gives \eqref{it:eq-two-transitive}~$\Rightarrow$~\eqref{it:eq-uniform-admissible}, and \eqref{it:eq-uniform-admissible}~$\Rightarrow$~\eqref{it:eq-admissible} is clear. The equivalence \eqref{it:eq-admissible}~$\Leftrightarrow$~\eqref{it:eq-ccr} is Lemma~\ref{lem:ccr-admissible}, and \eqref{it:eq-ccr}~$\Rightarrow$~\eqref{it:eq-type-i} follows from \cite{Ka51}. The implication \eqref{it:eq-type-i}~$\Rightarrow$~\eqref{it:eq-two-transitive} is \cite[Corollary~D]{CKM23}. The von Neumann algebra characterization of type~$\I$ groups gives \eqref{it:eq-type-i}~$\Rightarrow$~\eqref{it:eq-vn-type-i}, and every type~$\I$ von Neumann algebra is amenable, giving \eqref{it:eq-vn-type-i}~$\Rightarrow$~\eqref{it:eq-vn-amenable}.

    It remains to prove \eqref{it:eq-vn-amenable}~$\Rightarrow$~\eqref{it:eq-uniform-admissible}. Assume that $\LL(G)$ is amenable. Nonamenability of $G$ and Lemma~\ref{lem:tree-thick-reduction} allow us to replace $T$ by a thick locally finite tree $S$ on which $G$ acts minimally as a closed subgroup of $\Aut(S)$. By \cite[Theorem~C]{HR19}, this action is locally two-transitive. Put $H=G^+$. For adjacent vertices $v_1,v_2\in V(S)$, set $K_i=H_{v_i}=G_{v_i}$ for each $i\in\{1,2\}$ and $J=K_1\cap K_2$. Lemma~\ref{lem:tree-amalgam} gives
    $$
        H=K_1\ast_JK_2.
    $$

    Since $H$ is open in $G$, the inclusion $\LL(H)\subset\LL(G)$ is expected by Subsection~\ref{subsec:prelim-group-vna}, so $\LL(H)$ is amenable. Theorem~\ref{thm:compact-amalgam-ccr} shows that $H$ is uniformly admissible, and Fact~\ref{fact:uniform-finite-extension} gives the same conclusion for $G$.
\end{proof}

\appendix

\section*{Appendix. From \texorpdfstring{$\Wstar$}{W*}-amenability to boundary two-transitivity}
\refstepcounter{section}
\label{sec:amenability-boundary}

In this appendix, we prove the implication \eqref{it:eq-vn-amenable}~$\Rightarrow$~\eqref{it:eq-two-transitive} of Theorem~\ref{thm:equivalences}, which recovers \cite[Corollary~D]{CKM23}. The proof uses the local two-transitivity theorem of Houdayer and Raum \cite[Theorem~C]{HR19} and the Hecke algebra of a vertex stabilizer. The merger construction of Section~\ref{sec:atomic-afp} shows that a suitable spherical corner is diffuse abelian. Commutativity of the corresponding Hecke algebra then forces transitivity on spheres, and hence two-transitivity on the boundary.

\begin{The}
    \label{thm:amenability-boundary}
    Let $G$ be a closed nonamenable subgroup of $\Aut(T)$ acting minimally on a locally finite tree $T$. If $\LL(G)$ is amenable, then the action of $G$ on $\partial T$ is two-transitive.
\end{The}

We first isolate the operator algebraic consequence of having two double cosets in each compact factor. This condition yields a diffuse abelian spherical corner at the start of the merger construction, and amenability ensures that the subsequent mergers leave it unchanged.

\begin{Prop}
    \label{prop:diffuse-spherical-corner}
    Let $G=K_1\ast_JK_2$ be an amalgam of second countable compact groups over a common open subgroup $J$ and put $d_i=[K_i:J]$ for each $i\in\{1,2\}$. Assume that $d_1\geq d_2\geq2$, that $\abs{J\backslash K_i/J}=2$ for each $i\in\{1,2\}$, and that $\LL(G)$ is amenable. Then $p_{K_1}\LL(G)p_{K_1}$ is diffuse abelian.

    In particular, the Hecke algebra $\Hk G{K_1}$ is commutative.
\end{Prop}

\begin{proof}
    The group $G$ is unimodular by Subsection~\ref{subsec:prelim-bass-serre}. Normalize Haar measure by $\mu_G(J)=1$ and put
    $$
        M=\LL(G),\qquad A_i=\LL(K_i),\qquad B=\LL(J).
    $$
    By~\eqref{eq:compact-amalgam-afp}, we have $M=A_1\ast_BA_2$. Write $\tau=\tau_G$ and set
    $$
        p=p_J,\qquad e_i=p_{K_i},\qquad a_i=\tau(e_i)=\frac1{d_i}.
    $$
    Then $e_i\leq p$, $\tau(p)=1$ and $0<a_1\leq a_2\leq1/2$.

    We first identify the corners by $p$. Since $p$ is the averaging projection of $J$, it is a minimal central projection of $B$ and $pBp=\cc p$. The corner $pA_ip$ is the image under $\lambda_G$ of the finite-dimensional Hecke algebra $\Hk {K_i}J$. Its dimension is the number of double cosets of $J$ in $K_i$, hence is two. Since $e_i$ and $p-e_i$ are nonzero orthogonal projections in this corner,
    \begin{equation}\label{eq:spherical-local-corners}
        pA_ip=\cc e_i\oplus\cc(p-e_i).
    \end{equation}

    Choose $f\in B$ by taking one minimal projection in each matrix block of $B$, as in the proof of Theorem~\ref{thm:lattice}. The scalar block $\cc p$ contributes $p$, so $p\leq f$. Put $D=fBf$ and choose atomic masas $D_i\subset fA_if$ containing $D$. The full-corner identity gives
    $$
        fMf\cong fA_1f\ast_DfA_2f.
    $$
    Use the increasing sequence $(N_n)_{n\geq0}$ constructed in the proof of Proposition~\ref{prop:abelian-amalgam}, starting from $N_0=D_1\ast_DD_2$. To apply that construction, equip the algebras with the trace $\lcm(d_1,d_2)\tau$. Their atoms then have integral traces by Lemma~\ref{lem:compact-amalgam-afp}. This rescaling changes neither the conditional expectations nor the algebras. We retain the original trace $\tau$ for the computations below.

    Each $N_n$ is amenable and of type~$\I$, each nontrivial inclusion $N_n\subset N_{n+1}$ satisfies the hypotheses of Proposition~\ref{prop:amenable-free-merger} with the rescaled trace, and
    $$
        fMf=\Wstar\Bigl(\bigcup_{n\geq0}N_n\Bigr).
    $$
    Since $p\in D\subset D_i$, the algebra $pD_i$ is a masa in $pA_ip$. The latter is abelian by~\eqref{eq:spherical-local-corners}, so $pD_i=pA_ip$. In particular, $e_1,e_2\in N_0$. Also, $p$ is central in $N_0$, and compression by $p$ gives
    $$
        pN_0p\cong
        \bigl(\cc e_1\oplus\cc(p-e_1)\bigr)
        \ast\bigl(\cc e_2\oplus\cc(p-e_2)\bigr).
    $$

    Apply Proposition~\ref{prop:two-projections} with $\alpha=1-a_1$, $\beta=a_2$ and $\gamma=1$, ordering the first factor as $\cc(p-e_1)\oplus\cc e_1$. Its hypotheses hold because $a_1\leq a_2\leq1/2$. The diffuse summand of $pN_0p$ is
    $$
        \rM_2(\cc)\ovt\rL^\infty([0,1],\nu),
        \qquad\text{with trace }a_1(\Tr_2\pt\nu).
    $$
    On this summand, $e_1$ has rank one almost everywhere. Its component there therefore has trace $a_1=\tau(e_1)$. Faithfulness of $\tau$ shows that its atomic component is zero. Consequently, $e_1N_0e_1$ is diffuse abelian.

    For every $n\geq0$, successive applications of Claims~\ref{fact:merger-key-facts}\textup{(\ref*{it:merger-fact-diffuse})} give
    $$
        e_1N_ne_1=e_1N_0e_1.
    $$
    Taking the weak$^\star$ closure and using $e_1\leq f$, we obtain
    $$
        e_1Me_1=e_1(fMf)e_1=e_1N_0e_1.
    $$
    Thus $p_{K_1}\LL(G)p_{K_1}$ is diffuse abelian. The last assertion follows because $\lambda_G$ embeds $\Hk G{K_1}$ into this corner.
\end{proof}

Proposition~\ref{prop:ttn-spherical-amenability} starts from spherical transitivity and uses linear growth bounds for the Hecke algebras of smaller compact open subgroups to prove $\Wstar$-amenability. Here we recover spherical transitivity from commutativity of the Hecke algebra of a single vertex stabilizer, using the unique midpoint of a geodesic.

\begin{Lem}
    \label{lem:commutative-hecke-boundary}
    Let $T$ be a thick locally finite tree and let $G\leq\Aut(T)$ be a closed type-preserving subgroup. Fix $v\in V(T)$ and put $K=G_v$. Suppose that $G$ acts transitively on the bipartition class containing $v$, that is, on the vertices at even distance from $v$. If $\Hk GK$ is commutative, then $G$ acts two-transitively on $\partial T$.
\end{Lem}

\begin{proof}
    The subgroup $K$ is compact and open. Normalize Haar measure by $\mu_G(K)=1$ and put $X=Gv$. We view Hecke functions as matrices indexed by $X$. To each $f\in\Hk GK$ we associate the matrix $A_f$ defined, for all $g,h\in G$, by
    $$
        A_f(gv,hv)=f(g^{-1}h).
    $$
    This is well defined by bi-$K$-invariance. The matrix $A_f$ is invariant under the diagonal action of $G$ on $X\times X$. Indeed, for all $s,g,h\in G$,
    $$
        A_f(sgv,shv)=f((sg)^{-1}sh)=f(g^{-1}h)=A_f(gv,hv).
    $$
    Compact support of $f$ implies that the nonzero entries of $A_f$ occur only at pairs of vertices at bounded distance. Since $T$ is locally finite, every row and every column of $A_f$ is finitely supported. Conversely, every $G$-invariant matrix $A$ with finitely supported rows and columns gives a Hecke function $f$ defined by $f(g)=A(v,gv)$ for every $g\in G$. Decomposing the convolution integral into left cosets of $K$, each of measure one, gives, for all $f_1,f_2\in\Hk GK$,
    $$
        A_{f_1\ast f_2}=A_{f_1}A_{f_2}.
    $$
    Thus these matrices commute by hypothesis.

    We first prove that $K$ is transitive on every even sphere about $v$. Fix an even integer $\ell\geq2$ and vertices $x,y$ at distance $\ell$ from $v$. Both belong to $X$ by the orbit assumption. Let $A$ be the indicator matrix of pairs in $X\times X$ at distance $\ell$, and let $B$ be the indicator matrix of the orbit of $(v,x)$ under the diagonal action of $G$. Both matrices are $G$-invariant, and their rows and columns are supported on finite spheres, so $AB=BA$.

    Since $T$ is thick, choose a vertex $z$ at distance $\ell$ from $v$ in a branch different from those containing $x$ and $y$. Then $z\in X$, and both $[z,x]$ and $[z,y]$ have length $2\ell$ and midpoint $v$. For each $w\in\{x,y\}$, a nonzero term in either matrix product at $(z,w)$ requires an intermediate vertex at distance $\ell$ from both $z$ and $w$. The only such vertex is the midpoint $v$. Therefore
    $$
        \begin{aligned}
            (AB)(z,w)&=\sum_{u\in X}A(z,u)B(u,w)=B(v,w),\\
            (BA)(z,w)&=\sum_{u\in X}B(z,u)A(u,w)=B(z,v).
        \end{aligned}
    $$
    Commutativity gives
    $$
        B(v,y)=B(z,v)=B(v,x)=1.
    $$
    By the definition of $B$, some element of $G$ fixes $v$ and maps $x$ to $y$. This element belongs to $K$, proving transitivity on every even sphere.

    We now use compactness of $K$ to pass from spheres to the boundary. Fix $\xi,\eta\in\partial T$. For every $n\geq1$, let $x_n$ and $y_n$ be the vertices at distance $2n$ from $v$ on the rays to $\xi$ and $\eta$. The transitivity just proved shows that each set
    $$
        \{k\in K\mid kx_n=y_n\}
    $$
    is nonempty and closed in $K$. These sets are nested because an automorphism fixing $v$ and mapping $x_{n+1}$ to $y_{n+1}$ maps $[v,x_{n+1}]$ onto $[v,y_{n+1}]$, and hence maps $x_n$ to $y_n$. Compactness of $K$ gives an element in their intersection. This element maps every finite initial segment of the ray to $\xi$ onto the corresponding segment of the ray to $\eta$, so it maps $\xi$ to $\eta$. Thus $K$, and hence $G$, is transitive on $\partial T$.

    Finally, $Gv$ contains vertices at arbitrarily large even distances from $v$. Since a compact group has finite orbits on the discrete set $V(T)$, the group $G$ is noncompact. By \cite[Lemma~3.1.1]{BM00}, its action on $\partial T$ is therefore two-transitive.
\end{proof}

We can now combine these two results. Local two-transitivity supplies the double-coset condition in Proposition~\ref{prop:diffuse-spherical-corner}, while the amalgam decomposition supplies the vertex-orbit hypothesis in Lemma~\ref{lem:commutative-hecke-boundary}.

\begin{proof}[Proof of Theorem~\ref{thm:amenability-boundary}]
    Nonamenability of $G$ excludes cases~\textup{(\ref*{it:tree-finite})} and~\textup{(\ref*{it:tree-line})} of Lemma~\ref{lem:tree-thick-reduction}. By case~\textup{(\ref*{it:tree-thick})}, we may replace $T$ by a thick locally finite tree $S$ with a $G$-equivariant homeomorphism $\partial T\cong\partial S$, and identify $G$ with a closed subgroup of $\Aut(S)$. Since $\LL(G)$ is amenable, \cite[Theorem~C]{HR19} implies that the action on $S$ is locally two-transitive.

    For adjacent vertices $v_1,v_2\in V(S)$, let $K_i=(G^+)_{v_i}=G_{v_i}$ and $J=K_1\cap K_2$. By Lemma~\ref{lem:tree-amalgam},
    $$
        G^+=K_1\ast_JK_2,
    $$
    and $d_i=[K_i:J]$ is the degree of $v_i$. In particular, $d_i\geq3$. The orbit map $kJ\mapsto kv_2$ identifies $K_1/J$ with the neighbors of $v_1$, and the analogous statement holds for $K_2/J$. Local two-transitivity gives exactly two $J$-orbits on $K_i/J$, or equivalently exactly two double cosets of $J$ in $K_i$. Relabel the vertices so that $d_1\geq d_2$.

    Since $G^+$ is open in $G$, the inclusion $\LL(G^+)\subset\LL(G)$ is expected by Subsection~\ref{subsec:prelim-group-vna}. Hence $\LL(G^+)$ is amenable.

    Proposition~\ref{prop:diffuse-spherical-corner} shows that $\Hk{G^+}{K_1}$ is commutative. The group $G^+$ is closed and type-preserving, and its orbit $G^+v_1$ consists of all vertices at even distance from $v_1$. Lemma~\ref{lem:commutative-hecke-boundary} implies that $G^+$, and therefore $G$, acts two-transitively on $\partial S$, and hence on $\partial T$.
\end{proof}

\end{document}